\input amstex
\magnification = 1200

\define\Om{\Omega}
\define\om{\omega}

\define\ph{\varphi}
\define\rh{\rho}
\define\ps{\psi}
\define\Ph{\Phi}
\define\al{\alpha}
\define\be{\beta}
\define\Ga{\Gamma}
\define\ga{\gamma}
\define\de{\delta}
\define\ze{\zeta}
\define\De{\Delta}
\define\ep{\varepsilon}
\define\et{\eta}
\define\la{\lambda}
\define\ka{\varkappa}

\define\si{\sigma}
\define\ta{\tau}
\define\th{\theta}
\define\ch{\chi}

\define\dif{\bar\partial}
\define\ad{\operatorname{ad}}
\define\Ad{\operatorname{Ad}}
\define\sgn{\operatorname{\sgn}}
\define\id{\operatorname{id}}

\define\SL{\operatorname{SL}}
\define\GL{\operatorname{GL}}

\define\SO{\operatorname{SO}}
\redefine\Sp{\operatorname{Sp}}
\define\U{\operatorname U}

\define\Ker{\operatorname{Ker}}
\define\Int{\operatorname{Int}}
\define\Aut{\operatorname{Aut}}
\define\Hom{\operatorname{Hom}}
\define\Exp{\operatorname{Exp}}
\define\End{\operatorname{End}}

\define\Pic{\operatorname{Pic}}
\redefine\Im{\operatorname{Im}}

\documentstyle{amsppt}

\topmatter
\title
Triangular de Rham cohomology of compact K\"ahler manifolds
\endtitle
\rightheadtext{triangular de Rham cohomology}
\author
A. Brudnyi, A. Onishchik
\endauthor
\affil Ben Gurion University of the Negev, Yaroslavl State University
\endaffil
\address
P.O.B. 653, Beer-Sheva 84105, Israel
\endaddress
\email
brudnyi\@cs.bgu.ac.il
\endemail
\address
Yaroslavl State University, Sovetskaya 14, 150 000 Yaroslavl, Russia
\endaddress
\email
arkadiy\@onishchik.msk.ru
\endemail
\thanks
Work supported in part by the Russian Foundation for Fundamental Research
(Grant 98-01-00329) (the second author).
\endthanks
\abstract
We study the de Rham 1-cohomology $H^1_{DR}(M,G)$ of a smooth manifold $M$
with values in a Lie group $G$. By definition, this is the quotient of the
set of flat connections in the trivial principal bundle $M\times G$ by the
so-called gauge equivalence. We consider the case when $M$ is a compact
K\"ahler manifold and $G$ is a solvable complex linear algebraic group of a
special class which contains the Borel subgroups of all complex classical
groups and, in particular, the group $T_n(\Bbb C)$ of all triangular
matrices. In this case, we get a description of the set $H^1_{DR}(M,G)$ in
terms of the 1-cohomology of $M$ with values in the (abelian) sheaves of flat
sections of certain flat Lie algebra bundles with
fibre $\frak g$ (the Lie algebra of $G$) or, equivalently, in terms of the
harmonic forms on $M$ representing this cohomology.
\endabstract
\endtopmatter

\head
0. Introduction
\endhead

The paper is devoted to the study of the de Rham 1-cohomology $H^1_{DR}(M,G)$
of a smooth manifold $M$ with values in a Lie group $G$. If $G$ is
non-abelian, then $H^1_{DR}(M,G)$ admits no natural group structure and is
usually regarded as a set with a distinguished point. By definition, this is
the quotient of the set of flat connections in the trivial principal bundle
$M\times G$ by the so-called gauge equivalence, the distinguished point being
the class of the zero connection. Note that the de Rham 1-cohomology set has
two important interpretations. First, $H^1_{DR}(M,G)$ admits a natural
injection into the \v Cech 1-cohomology set $H^1(M,G)$, the image being
interpreted as the set of smoothly (or topologically) trivial flat
principal bundles with base $M$ and structure group $G$. In the case $G =
\Bbb R$, this correspondence is the classical de Rham isomorphism. Second,
any flat connection determines the holonomy homomorphism $\pi_1(M)\to G$,
giving rise to an injective mapping
$H^1_{DR}(M,G)\to\Hom(\pi_1(M),G)/\Int G$. If $G$ is contractible (e.g.,
solvable and simply connected), then any smooth principal bundle with
structure group $G$ is trivial, and both injections are bijections.

We consider the case when $M$ is a compact K\"ahler manifold and $G$ is a
solvable complex linear algebraic group of a special class which contains
the Borel subgroups of all complex classical groups and, in particular, the
group $T_n(\Bbb C)$ of all non-singular triangular $n\times n$-matrices. In
this case, we get a description of the set
$H^1_{DR}(M,G)$ in terms of the 1-cohomology of $M$ with values in the
(abelian) sheaves of flat sections of certain flat Lie algebra bundles with
fibre $\frak g$ (the Lie algebra of $G$) or, equivalently, in terms of the
harmonic forms on $M$ representing this cohomology. Our method is based on
using the non-abelian Dolbeault cohomology set $H^{0,1}(M,G)$ of $M$ with
values in $G$ that is defined similarly to the de Rham cohomology set and
essentially exploits the properties of the Hodge decomposition on a compact
K\"ahler manifold.

The paper is devided into two sections. The first one contains the necessary
facts on non-abelian de Rham and Dolbeault cohomology. In particular, the
technics of twisting described here is very important for the sequel. The
proofs are omitted, referring the reader to [O3, O4]. We also expose some
facts on Hodge theory for flat vector bundles.

Section 2 contains formulations and proofs of the main results of the paper.
We introduce the so-called Hodge property for solvable complex algebraic
groups. Lemma 2.2 gives a method to construct groups having this property,
while Theorem 2.1 presents a list of such groups including the Borel
subgroups of classical complex linear groups (one of these subgroups is
$T_n(\Bbb C)$). The Hodge
property for $G$ implies, in particular, that the natural mapping
$\Pi_{0,1}^*: H^1_{DR}(M,G)\to H^{0,1}(M,G)$ is surjective whenever $M$ is
compact and K\"ahler. Theorem 2.2 describes the so-called canonical
representatives of de Rham cohomology classes in terms of harmonic forms with
values in certain Lie algebra bundles. The final result about $H^1_{DR}(M,G)$
is formulated as Theorem 2.3. In the last subsection, we formulate two
results that can be proved by the same argument as Theorem 2.3. Theorem 2.4
gives a description of a twisted version of the de Rham cohomology of a
compact K\"ahler manifold $M$ with values in the unipotent radical of a
group $G$ having the Hodge property. In Theorem 2.5, the situation is
considered when such a group $G$ is a subgroup of an algebraic group
$\hat G$. Here we describe the quotient of the set of flat connections $\om$
on $M\times\hat G$ such that $\Pi_{0,1}\om$ takes its values in $\frak g$ by
the gauge equivalence determined by $G$.

In the case when $G = T_n(\Bbb C)$, the main results of the paper were
proved in the research thesis [Br1] of the first author. The thesis
also contains applications of these results to the study of the fundamental
group $\pi_1(M)$ of a compact K\"ahler manifold $M$ and, in particular, a
classification of compact solvmanifolds that admit K\"ahler structures. These
applications will be published in [Br2].

\head
1. De Rham and Dolbeault cohomology with values in a Lie group
\endhead

{\bf 1.1.}
Here we discuss a non-linear cochain complex that coincides, in the classical
abelian case, with the initial part of the usual de Rham complex of a smooth
manifold. We follow [O3], Sections 4 and 5, and [O4], Section 2.

Let us first introduce some notation. Let $X$ be a topological space,
$\Cal F$ a sheaf of (non-necessary abelian)
groups on $X$. Then one defines the 0-cohomology group $H^0(X,\Cal F)$ which
coincides with the group of global sections $\Ga(X,\Cal F)$. One also defines
the \v Cech 1-cohomology $H^1(X,\Cal F)$ which in general is not a group, but
merely a set with a distinguished point. Let $\frak U$ be an open cover of
$X$. In a usual way, one defines the groups of $p$-cochains
$C^p(\frak U,\Cal F),\; p\ge 0$, and the set of 1-cocycles
$$
Z^1(\frak U,\Cal F) = \{z\in C^1(\frak U,\Cal F)\mid z_{ij}z_{jk} = z_{ik}
\text{ in } U_i\cap U_j\cap U_j\ne\emptyset\}.
$$
There is an action $\rh$ of $C^0(\frak U,\Cal F)$ on $C^1(\frak U,\Cal F)$
given by
$$
(\rh(a)(c))_{ij} = c_iz_{ij}c_j^{-1}.
$$
The set $Z^1(\frak U,\Cal F)$ is invariant under $\rh$. Forming the quotient
$H^1(\frak U,\Cal F) = \mathbreak
Z^1(\frak U,\Cal F)/\rh(C^0(\frak U,\Cal F))$ and
passing to a limit over all open covers $\frak U$, one obtains the cohomology
set $H^1(X,\Cal F)$. Its distinguished point $\ep$ is the class of the unit
cocycle $e\in Z^1(\frak U,\Cal F)$.

Let now $G$ be a Lie group. For any smooth manifold $M$, denote by
$\Cal F_G$ the sheaf of germs of smooth $G$-valued functions on $M$. In
particular, $\Cal F = \Cal F_{\Bbb R}$ is the structure sheaf of $M$.
Clearly, $\Cal F_G$ is a sheaf of groups. The 0-cohomology $H^0(M,\Cal F_G)$
is the group $F_G$ of global smooth functions $M\to G$, while the
1-cohomology set $H^1(M,\Cal F_G)$ is usually identified with the set of all
smooth principal bundles with base $M$ and structure group $G$ (regarded up
to isomorphy). Namely, if $\frak U = (U_i)$ is an open cover of $M$, then any
$z\in Z^1(\frak U,\Cal F_G)$ determines the principal bundle with the
transition functions $z_{ij}$ got from $M\times G$ by twisting with the help
of $z$; this is the bundle corresponding to the cohomology class $\ze\in
H^1(M,\Cal F_G)$ of the cocycle $z$. The unit element $\ep\in
H^1(M,\Cal F_G)$ corresponds to the trivial bundle $M\times G$.

Clearly, we may identify the constant sheaf $G$ on $M$ with the subsheaf of
$\Cal F_G$ consisting of germs of {\it flat} (i.e., locally constant)
functions. Now, $H^0(M,G)$ is the group of all flat functions $M\to G$. The
set $H^1(M,G)$ will be interpreted as the set of {\it flat} principal bundles
with base $M$ and structure group $G$, i.e., principal bundles with locally
constant transition functions, regarded up to corresponding isomorphy.

Let $\frak g$ be the tangent Lie algebra of $G$ and
$\Ph_{\frak g} = \bigoplus_{p\ge 0}\Ph^p_{\frak g}$ the graded sheaf of
$\frak g$-valued smooth forms on $M$. This is a sheaf of graded Lie
superalgebras, the bracket $[\,,\,]$ being induced by that of $\frak g$.
Denoting $A^p_{\frak g} = \Ga(G,\Ph^p_{\frak g})$, we obtain the usual
$\frak g$-valued de Rham complex $(A_{\frak g},d)$, where $A_{\frak g} =
\bigoplus_{p\ge 0}A^p_{\frak g}$ has a structure of graded Lie superalgebra
and the exterior derivative $d$ is its derivation of degree 1.

In order to define the desired non-linear complex, we denote by $\varpi\in
A^1_{\frak g}$ the canonical 1-form on $G$ assigning to any tangent
vector $v$ at $g\in G$ the vector $d_gr_g^{-1}(v)\in\frak g$, where $r_g:
x\mapsto xg$ is the right translation corresponding to $g$. Consider the
triple of groups $R_G = \{R^0_G,R^1_G,R^2_G\} = \{F_G,A^1_{\frak g},
A^2_{\frak g}\}$ and define the coboundary operators $\de_0: R^0_G\to R^1_G$
and $\de_1: R^1_G\to R^2_G$ by
$$
\align
\de_0(g) &= g^*(\varpi),\ \  g\in F_G,\tag1.1 \\
\de_1(\al) &= d\al -\frac12[\al,\al],\ \ \al\in A^1_{\frak g}.\tag1.2
\endalign
$$
Note that $R^0_G = F_G$ acts on the graded algebra $A_{\frak g}$ by
automorphisms induced by the adjoint representation $\Ad$ of $G$. Denoting
this action by $\Ad$, we have
$$
d(\Ad g(\al)) = \Ad g(d\al) + [\de_0(g),\Ad g(\al)].\tag1.3
$$
One checks easily that $\de_0: F_G\to A^1_{\frak g}$ is a {\it crossed
homomorphism} with respect to $\Ad$, i.e.,
$$
\de_0(gh) = \de_0(g) + \Ad g(\de_0(h)).\tag1.4
$$
This gives rise to the following affine action of $R^0_G$ on $R^1_G =
A^1_{\frak g}$ by so-called {\it gauge transformations}:
$$
\rh(g)(\al) =  \Ad g(\al) + \de_0(g),\tag1.5
$$

In what follows, $G$ often will be a Lie subgroup of the linear group
$\GL_m(\Bbb R)$. In this case, the coboundary operators $\de_p$ and the
actions $\Ad$ and $\rh$ can be written as follows:
$$
\align
\de_0(g) &= (dg)g^{-1},\\
\de_1(\al) &= d\al -\al\wedge\al,\\
\Ad g(\al) &= g\al g^{-1},\\
\rh(g)(\al) &=  g\al g^{-1} + (dg)g^{-1}.
\endalign
$$

Using (1.3), one easily verifies the following relation:
$$
\de_1\circ\rh(g) = \Ad g\circ\de_1.\tag1.6
$$
This means that the triple $R_G$ together with the mappings $\de_p,\;
p = 0,1$, and the actions $\Ad$ and $\rh$ is a {\it non-abelian cochain
complex} in the sense of [O3, O4]. It is called the {\it de Rham complex
with values in} $G$.

Introduce the sets of cocycles
$$
Z^p(R_G) = \Ker\de_p = \de_p^{-1}(0),\ \ p = 0,1.
$$
Then $Z^0(R_G)$ is a subgroup of $F_G$. Now, the subset $Z^1(R_G)$ is
invariant under $\rh$ due to (1.6). We define the
{\it de Rham cohomology of} $M$ {\it with values in} $G$ (or the cohomology
of the complex $R_G$) by
$$
\align
H^0_{DR}(M,G) &= H^0(R_G) = Z^0(R_G),\\
H^1_{DR}(M,G) &= H^1(R_G) = Z^1(R_G)/\rh(F_G).
\endalign
$$
Then $H^0_{DR}(M,G) = H^0(M,G)$. If $G$ is non-abelian, then the set
$H^1_{DR}(M,G)$ does not admit any natural group structure. We regard it as a set with the
distinguished point $\ep = \rh(F_G)(0)$. The cohomology class of a cocycle
$\om\in Z^1(R_G)$ will be denoted by $[\om]$.

Note that on the sheaf level we get the following exact sequence of sheaves:
$$
e\to G\overset i\to\to\Cal F_G\overset\de_0\to\to\Ph^1_{\frak g}
\overset\de_1\to\to\Ph^2_{\frak g},\tag1.7
$$
where $i$ is the natural injection. It implies the following relation
between the de Rham and the \v Cech 1-cohomology (see [O3], Section 5).

\proclaim{Proposition 1.1}
We have the following exact sequence of sets with distinguished points:
$$
e\to H^1_{DR}(M,G)\overset\mu\to\to H^1(M,G)\overset i_1^*\to\to
H^1(M,\Cal F_G).\tag1.8
$$
Here $i_1^*$ is determined by $i$, while $\mu$ is defined as follows. For any
cocycle $\om\in Z^1(R_G)$,
choose an open cover $\frak U = (U_i)$ of $M$ such that $\om = \de_0(c_i)$ in
any $U_i$ for certain smooth $c_i: U_i\to G$. Then $z_{ij} = c_i^{-1}c_j$
form a cocycle $z = (z_{ij})\in Z^1(\frak U,G)$, and $\mu$ sends $[\om]$ to
the cohomology class of $z$.
\endproclaim

The exactitude of (1.8) means that $\mu$ maps bijectively the set
$H^1_{DR}(M,G)$ onto the subset of
those flat bundles from $H^1(M,G)$ which are trivial as smooth bundles.

In the case $G = \Bbb R$, the complex $R_G$ is, clearly, a part of the
classical de Rham complex, and (1.8) gives the de Rham theorem for the
1-cohomology.

If $G$ is non-abelian, then (1.8) is not an exact sequence of groups and
their homomorphisms, and so the fibres of the mapping $i_1^*$ cannot be
described in terms of its kernel. To describe these fibres, it is necessary
to consider twisted versions of the de Rham complex.

Now we discuss certain twisting constructions. Let $\Cal S$ be
a sheaf of sets on $X$ and let us denote by $\Cal Aut\,\Cal S$ the sheaf of
germs of automorphisms of this sheaf. Its sections over an open set $U\subset
X$ are, by definition, bijective continuous mappings $\Cal S|U\to\Cal S|U$
leaving each stalk invariant. Clearly, $\Cal Aut\,\Cal S$ is a
sheaf of groups.
Let $\frak U = (U_i)$ be an open cover of $X$, and fix a cocycle $z\in
Z^1(\frak U,\Cal Aut\,\Cal S)$. Then we can twist $\Cal S$ with the help of
$z$ getting a new sheaf $\Cal S^z$ on $X$. This means that we glue together
any two sheaves $\Cal S|U_i$ and $\Cal S|U_j$ over $U_i\cap U_j\ne\emptyset$
identifying $a_i\in (\Cal S|U_i)_x$ and $a_j\in (\Cal S|U_j)_x,\;
x\in U_i\cap U_j$, under the condition $a_i = z_{ij}(a_j)$. If $S$ is a sheaf
of groups and we use automorphisms of a sheaf of groups, then, clearly,
$\Cal S^z$ is a sheaf of groups, too, etc.

The {\it twisted sheaf} $\Cal S^z$ depends, up to isomorphism, only on the
cohomology class $\ze\in H^1(X,\Cal Aut\,\Cal S)$ determined by $z$. E.g., if
$z$ is cohomologous to $e$, i.e., $z_{ij} = c_ic_j^{-1}$ for a 0-cochain $c
\in C^0(\frak U,\Cal Aut\,\Cal S)$, then $a_i = z_{ij}(a_j)$ implies
$c_i^{-1}(a_i) = c_j^{-1}(a_j)$. Thus, the correspondence $(a_i)\mapsto
(c_i^{-1}(a_i))$ is an isomorphism of $\Cal S^z$ onto $\Cal S$.

Note that any section $s\in\Ga(U,\Cal S^z)$ is given by a 0-cochain
$(s_i)$, where $s_i\in\Ga(U\cap U_i,\Cal S)$, satisfying $s_i = z_{ij}(s_j)$
over $U\cap U_i\cap U_j$.

To twist the de Rham complex with values in a Lie group $G$, suppose that a
cocycle $\frak z\in Z^1(\frak U,\Aut G)$ is given, $\Aut G$ being the group
of automorphisms of $G$. Consider the twisted sheaf of groups $G^{\frak z}$.
We can realize it as the sheaf $\Cal C_E$ of locally constant sections of the
flat group bundle $E$ which is got from $M\times G$ by twisting with the help
of $\frak z$.
Any automorphism of $G$ induces an automorphism of $\Cal F_G$, and hence we
get the twisted sheaf $(\Cal F_G)^{\frak z}$ which is the sheaf $\Cal F_E$ of
smooth sections of $E$. We also have the cocycle $d\frak z$ of automorphisms
of the Lie algebra $\frak g$ giving rise to the twisted sheaves
$(\Ph^p_{\frak g})^{d\frak z},\; p\ge 0$; these are the sheaves
$\Ph^p_{\frak e}$ of $p$-forms with values in the Lie algebra bundle
$\frak e$ which is got from $M\times\frak g$ by twisting with the help of
$d\frak z$.

Define now the graded Lie superalgebra $A_{\frak e}=
\bigoplus_{p\ge 0}A^p_{\frak e}$, where $A^p_{\frak e} =
\Ga(M,\Ph^p_{\frak e})$. Setting
$$
(d\al)_i = d\al_i,
$$
we correctly define a derivation $d$ of $A_{\frak e}$ giving rise to the
twisted de Rham complex $(A_{\frak e},d)$. Its cohomology is denoted by
$H^p_{DR}(M,\frak e)$, by a generalized de Rham theorem it is isomorphic to
$H^p(M,\Cal C_{\frak e})$. The non-linear twisted complex is defined as the
triple $R_E = \{R^0_E,R^1_E,R^2_E\}$, where $R^0_E = F_E = \Ga(M,\Cal F_E),\;
R^p_E = A^p_{\frak e},\; p = 1,2$. The coboundary operators and the actions
of $F_E$ are correctly defined by
$$
\aligned
(\de_0(f))_i &= (\de_0(f_i))\text{ for }f = (f_i)\in F_E,\\
(\de_1(\al))_i &= (\de_1(\al_i))\text{ for }\al = (\al_i)\in A^1_{\frak e},\\
(\Ad f(\al))_i&=(\Ad f_i(\al)_i)\text{ for }\al = (\al_i)\in A^p_{\frak e},\\
(\rh(f)(\al))_i&=(\rh(f_i)(\al)_i)\text{ for }\al = (\al_i)\in A^1_{\frak e}.
\endaligned\tag1.9
$$
As above, we obtain a non-abelian cochain complex. One defines the
0-cohomology group $H^0_{DR}(M,E) = H^0(R_E) = Z^0(R_E)$ and the 1-cohomology set $H^1_{DR}(M,E) = H^1(R_E) =
Z^1(R_E)/\rh(F_E)$. The group $H^0_{DR}(M,E)$ coincides with the group
\linebreak $H^0(M,\Cal C_E)$ of flat sections of $E$. Generalizing (1.8), we
get the following exact sequence of sets with distinguished points:
$$
e\to H^1_{DR}(M,E)\overset\mu\to\to H^1(M,\Cal C_E)\overset i_1^*\to
\to H^1(M,\Cal F_G).\tag1.10
$$

Let us denote by $\Int h$ the inner automorphism $x\mapsto hxh^{-1}$ of a
group $H$ determined by $h\in H$.
Suppose that $\frak z = \Int z$, i.e., $\frak z_{ij} = \Int {z_{ij}}$,
where $z = (z_{ij})\in Z^1(\frak U,G)$. Denote by $P$ the flat principal
bundle with the structure group $G$ over $M$ corresponding to $z$. Then $E =
\Int P$ is the group bundle associated to $P$ by the action $\Int$ of $G$ on
itself by inner automorphisms, and $\frak e = \Ad P$ is the Lie algebra
bundle associated to $P$ by the adjoint representation $\Ad$ of $G$ on
$\frak g$. The vector space $R^1_{\Int P} = A^1_{\Ad P}$
can be interpreted as the space of connections in $P$ (regarded as a smooth
principal bundle), while $\de_1(\al)\in R^2_{\Int P} = A^2_{\Ad P}$ is
interpreted as the curvature of a connection $\al$, so that $Z^1(R_{\Int P})$
is the set of connections with zero curvature (or {\it flat} connections).
Now, $R^0_{\Int P} = F_{\Int P}$ is identified with the group of smooth
automorphisms of the principal bundle $P$ (inducing the identity mapping of
the base $M$), and $\rh$ is the action of this group on $A^1_{\Ad P}$ by
gauge transformations. By definition, $H^1(R_{\Int P}) = H^1_{DR}(M,\Int P)$
is the quotient of $Z^1(R_{\Int P})$ by this action. Also $H^0(R_{\Int P}) =
H^0(M,\Cal C_{\Int P})$ is the group of flat automorphisms of $P$.
There exists a natural bijection of $H^1(R_{\Int P})$ onto the
subset $(i_1^*)^{-1}(i_1^*(\ze))\subset H^1(M,G)$ of all flat principal
bundles isomorphic to $P$ as smooth bundles (or, which is the same, onto the
set of all flat structures in the bundle $P$ regarded as a smooth principal
bundle, up to automorphisms of this smooth bundle), see [O3], Section 5.

Now we consider the special case when $P$ is trivial (as a smooth bundle).
This means that the cohomology class $\ze$ of $z$ lies in $\Im\mu$ (see
(1.8)). More precisely, we may suppose that $z_{ij} = c_i^{-1}c_j$, where
$c_i\in\Ga(U_i,\Cal F_G)$. Then $\ga = \de_0(c_i)$ is a well-defined global
form from $Z^1(R_G)$ such that $[\ga]$ satisfies $\mu([\ga]) = \ze$.
The cochain $(c_i)$ determines an isomorphism of groups $t:
R^0_{\Int P}\to R^0_G = F_G$ given by
$$
t((g_i)) = c_ig_ic_i^{-1}\ \ \text{in}\ \ U_i.\tag1.11
$$
Similarly, we get an isomorphism of graded Lie superalgebras $\ta:
A_{\Ad P}\to A_{\frak g}$ given on the $p$-components by
$$
\ta_p((\al_i)) = \Ad c_i(\al_i)\ \ \text{in}\ \ U_i,\ \ (\al_i)\in
A^p_{\Ad P}.\tag1.12
$$
This is not an isomorphism of complexes. More precisely, we see from (1.3)
that
$$
\ta\circ d = (d - \ad\ga)\circ\ta,
$$
where $\ad\ga: \al\mapsto [\ga,\al]$ is the adjoint operator in
$A_{\frak g}$ associated with $\ga$. Thus, the following is true:

\proclaim{Proposition 1.2}
The mapping $\ta$ defined by (1.12) is an isomorphism of complexes
$(A_{\Ad P},d)\to (A_{\frak g},d - \ad\ga)$.
\endproclaim

Now we establish a relation between the 1-cohomology sets of $R_G$ and
$R_{\Int P}$ (in the case when $P$ is trivial as a smooth bundle). Define
$r_{\ga}: A^1_{\Int P}\to A^1_{\frak g}$ by
$$
r_{\ga}((\al_i)) = \ta_1(\al_i) + \ga = \Ad c_i(\al_i) + \ga =
\rh(c_i)(\al_i).\tag1.13
$$
Clearly, this is a well defined bijective affine mapping.

\proclaim{Proposition 1.3}
In the notation given by (1.11), (1.12), (1.13), the following relations
hold:
$$
\align
\de_1&\circ r_{\ga} = \ta_2\circ\de_1,\\
r_{\ga}&\circ\rh(g) = \rh(t(g))\circ r_{\ga},\ \ g\in F_{\Int P}.
\endalign
$$
Thus, $r_{\ga}$ maps $Z^1(R_{\Int P})$ onto $Z^1(R_G)$ and induces a
bijection $r_{\ga}^*: H^1_{DR}(M,\Int P)\mathbreak\to H^1_{DR}(M,G)$ taking
$\ep$ to $[\ga]$.
\endproclaim

The proof is straightforward. The statement is also implied by Proposition
1.5 of [O4].

The construction of the de Rham complex gives a covariant functor $G\mapsto
R_G$ from the category of Lie groups into that of cochain complexes (see
[O4]). In fact,
to any smooth homomorphism of Lie groups $f: G\to Q$ there corresponds the
triple $\{f_0,f_1,f_2\}$ of homomorphisms $f_p: R^p_G\to R^p_Q$, where
$f_0(g) = f\circ g,\; g\in R^0_G = F_G$, and $f_p(\al) = df\circ\al,\;\al\in
R^p_G = A^p_{\frak g},\;p = 1,2$, satisfying
$$
\align
\de_p\circ f_p &= f_{p+1}\circ\de_p,\; p = 0,1,\\
f_1\circ\rh(g) &= \rh(f_0(g)) f_1.
\endalign
$$
Clearly, $f$ determines a homomorphism of groups $f_0^*: H^0_{DR}(M,G)\to
H^0_{DR}(M,Q)$ and a homomorphism of sets with distinguished points $f_1^*:
H^1_{DR}(M,G)\to\mathbreak H^1_{DR}(M,Q)$.

\medskip

{\bf 1.2.} Here we consider a theory, analogous to that exposed in
n$^{\circ} 1.1$, in which
flat bundles and locally constant sections or functions are replaced by
holomorphic ones. The abelian model is the classical Dolbeault complex of a
complex analytic manifold. We follow [O3], Section 6.

Let $G$ be a complex Lie group. For any complex manifold $M$, denote by
$\Cal O_G$ the sheaf of germs of $G$-valued holomorphic functions on $M$. In
particular, $\Cal O = \Cal O_{\Bbb C}$ is the structure sheaf of $M$.
Clearly, $H^0(M,\Cal O_G)$ is the group of all $G$-valued holomorphic
functions on $M$, while $H^1(M,\Cal O_G)$ can be interpreted as the set of
all principal holomorphic bundles with the base $M$ and the structure group
$G$, considered up to holomorphic isomorphisms leaving any point of $M$
fixed. The unit cohomology class $\ep$ corresponds to the trivial bundle
$M\times G$. As in Section 1.1, we will denote by $\Cal F_G$ the sheaf of
smooth $G$-valued functions on $M$ and by $\Ph_{\frak g}$ the sheaf of
$\frak g$-valued smooth forms on $M$. This
is a sheaf of bigraded complex Lie superalgebras, the bigrading being defined
by $\Ph_{\frak g} = \bigoplus_{p,q\ge 0}\Ph_{\frak g}^{p,q}$, where
$\Ph_{\frak g}^{p,q}$ is the sheaf of forms of type $(p,q)$. Let $\Pi_{p,q}:
\Ph^{p+q}\to\Ph^{p,q}$ be the natural projection. As usually, denote
$$
\partial\al = \Pi_{p+1,q}d\al,\;\dif\al =\Pi_{p,q+1}d\al,\ \ \al\in\Ph_{p,q},
$$
Then
$$
d = \partial + \dif,\ \ \partial^2 = \dif^2 = [\partial,\dif] = 0.
$$

Denote $A^{p,q}_{\frak g} = \Ga(M,\Ph_{\frak g}^{p,q})$. Then $A_{\frak g} =
\bigoplus_{p,q\ge 0}A^{p,q}_{\frak g}$ is a bigraded complex Lie
superalgebra, and $\partial$ and $\dif$ are derivations of bidegree $(1,0)$
and $(0,1)$, respectively. We get the bigraded Dolbeault complexes
$(A_{\frak g},\dif)$ and $(A_{\frak g},\partial)$. Denote by
$$
\align
\Om^p_{\frak g} &= \{\al\in\Ph_{\frak g}^{p,0}\mid\dif\al = 0\},\\
\overline\Om^p_{\frak g} &= \{\al\in\Ph_{\frak g}^{0,p}\mid\partial\al = 0\}
\endalign
$$
the subsheaves of holomorphic and antiholomorphic $p$-forms, respectively.
We have the classical Dolbeault isomorphisms
$$
\align
H^{p,q}(M,\frak g)\overset\operatorname{def}\to = H^q(A_{\frak g}^{p,*},\dif)
&\simeq H^q(M,\Om^p_{\frak g}),\tag 1.14 \\
\overline H^{p,q}(M,\frak g)\overset\operatorname{def}\to =
H^p(A_{\frak g}^{*,q},\partial)
&\simeq H^p(M,\overline\Om^q_{\frak g}).\tag1.15
\endalign
$$

Now we define the non-linear Dolbeault complex with values in $G$. Consider
the triple of groups $\overline R_G =
\{\overline R^0_G,\overline R^1_G,\overline R^2_G\} =
\{F_G,A_{\frak g}^{0,1},A_{\frak g}^{0,2}\}$ and define the coboundary
operators $\bar\de_0: \overline R^0_G\to\overline R^1_G$, $\bar\de_1:
\overline R^1_G\to\overline R^2_G$ and the action $\bar\rh$ of
$\overline R^0_G$ on $\overline R^1_G$ by
$$
\align
\bar\de_0(g) &= \Pi_{0,1}g^*(\varpi),\\
\bar\de_1(\al) &= \dif\al -\frac12[\al,\al],\\
\bar\rh(g)(\al) &= \Ad g(\al) + \bar\de_0(g).
\endalign
$$
Then we get a non-abelian cochain complex in the sense of [O3, O4], called
the {\it Dolbeault complex with values in} $G$. If $G$ is a Lie subgroup of
$\GL_m(\Bbb C)$, then
$$
\align
\bar\de_0(g) &= (\dif g)g^{-1},\\
\bar\de_1(\al) &= \dif\al - \al\wedge\al,\\
\bar\rh(g)(\al) &=  g\al g^{-1} + (\dif g)g^{-1}.
\endalign
$$
The cocycles of this complex are defined by
$$
Z^p(\overline R_G) = \Ker\bar\de_p = \bar\de_p^{-1}(0),\ \ p = 0,1.
$$
Then $Z^0(\overline R_G)$ is a subgroup of $F_G$. Now, the subset
$Z^1(\overline R_G)$ is invariant under $\bar\rh$. We define the
{\it Dolbeault cohomology of} $M$ {\it with values in} $G$ (or the cohomology
of $\overline R_G$) by
$$
\align
H^0(\overline R_G) &= Z^0(\overline R_G),\\
H^{0,1}(M,G) &= H^1(\overline R_G) = Z^1(\overline R_G)/\bar\rh(F_G).
\endalign
$$

Note that on the sheaf level we get the following exact sequence of sheaves:
$$
e\to\Cal O_G\overset i\to\to\Cal F_G\overset\bar\de_0\to\to
\Ph^{0,1}_{\frak g}\overset\bar\de_1\to\to\Ph^{0,2}_{\frak g},
$$
where $i$ is the natural injection. It implies the following relation
between the Dolbeault and the \v Cech 1-cohomology (see [O3], Section 6).

\proclaim{Proposition 1.4}
We have the following exact sequence of sets with distinguished points:
$$
e\to H^{0,1}(M,G)\overset\bar\mu\to\to H^1(M,\Cal O_G)\overset i_1^*\to\to
H^1(M,\Cal F_G).
$$
Here $i_1^*$ is determined by $i$, while $\bar\mu$ is defined as follows. For
any
cocycle $\om\in Z^1(\overline R_G)$, choose an open cover $\frak U = (U_i)$
of $M$ such that $\om = \bar\de_0(c_i)$ in any $U_i$ for certain smooth $c_i:
U_i\to G$. Then $z_{ij} = c_i^{-1}c_j$ determine a cocycle $z = (z_{ij})\in
Z^1(\frak U,\Cal O_G)$, and $\bar\mu$ takes the Dolbeault cohomology class
$[\om]$ of $\om$ to the cohomology class of $z$.
\endproclaim

This means that $\bar\mu$ maps bijectively the set $H^{0,1}(M,G)$ onto the
subset of those holomorphic bundles from $H^1(M,\Cal O_G)$ which are trivial
as smooth bundles.

As in Section 1.1, there exist twisted versions of the Dolbeault complex.
Denote by $\Aut_h G$ the group of holomorphic automorphisms of $G$ and choose
a cocycle $\frak z\in Z^1(\frak U,\Cal F_{\Aut_h G})$. Clearly,
$\Cal O_{\Aut_h G}$ acts on $\Cal O_G$ and $\Cal F_G$. We get the
holomorphic group bundle $E$ obtained by twisting $M\times G$ with the help
of $\frak z$ and the sheaves $\Cal O_E = (\Cal O_G)^{\frak z}$ and
$\Cal F_E = (\Cal F_G)^{\frak z}$ of holomorphic and smooth
sections of $E$, respectively. The cocycle $d\frak z$ gives rise to the
twisted sheaves $(\Ph^{p,q}_{\frak g})^{d\frak z},\; p,q\ge 0$; these are the
sheaves $\Ph^{p,q}_{\frak e}$ of $(p,q)$-forms with values in the holomorphic
Lie algebra bundle $\frak e$ which is got from $M\times\frak g$ by twisting
with the help of $d\frak z$. Clearly, the operator $\dif$ is well-defined in
$\Ph_{\frak e} = \bigoplus_{p,q\ge 0}\Ph^{p,q}_{\frak e}$, and one can
define the subsheaf of holomorphic $p$-forms $\Om^p_{\frak e}\subset
\Ph^{p,0}_{\frak e}$. Now, one can consider the bigraded Dolbeault complex
$(A_{\frak g},\dif)$, where $A_{\frak g} = \Ga(M,\Ph_{\frak e})$. Then
there are the Dolbeault isomorphisms similar to (1.14):
$$
H^{p,q}(M,\frak e)\overset\operatorname{def}\to =
H^q(A_{\frak e}^{p,*},\dif)\simeq H^q(M,\Om^p_{\frak e}).\tag 1.16
$$

The twisted version of the Dolbeault complex $\overline R_G$ is
$\overline R_E = \{\overline R^0_E,\overline R^1_E,\overline R^2_E\}$, where
$\overline R^0_E = F_E,\; \overline R^p_E =
A^{0,p}_{\frak e},\; p = 1,2$. The coboundary operators and the actions of
$F_E$ are correctly defined by formulas similar to (1.9). The 0-cohomology
of this cochain complex is the group $\Ga(M,\Cal O_E)$ of all holomorphic
sections of $E$; the 1-cohomology will be denoted $H^{0,1}(M,E)$. Let us
consider the case when $\frak z = \Int z$, where $z = (z_{ij})\in
Z^1(\frak U,\Cal O_G)$. Denote by $P$ the principal holomorphic bundle over
$M$ corresponding to $z$. Then $E = \Int P$ is the group bundle associated to
$P$ by the action $\Int$ of $G$ on itself by inner automorphisms, and
$\frak e = \Ad P$ is the Lie algebra bundle associated to $P$ by the adjoint
representation $\Ad$ of $G$ on $\frak g$.

Now we consider the special case when $P$ is trivial as a smooth bundle. We
may suppose that $z_{ij} = c_i^{-1}c_j$, where $c_i\in\Ga(U_i,\Cal F_G)$.
Then $\ga = \bar\de_0(c_i)$ is a well-defined global form from
$Z^1(\bar R_G)$, and its cohomology class $[\ga]$ satisfies $\bar\mu([\ga])
= \ze$, where $\ze$ is the cohomology class of $z$. Precisely as in Section
1.1, we get an isomorphism of groups $t: F_{\Int P}\to F_G$ given by (1.11)
and an isomorphism of bigraded Lie
superalgebras $\ta: A_{\Ad P}\to A_{\frak g}$ given by (1.12). Define a
mapping $\bar r_{\ga}: A^{0,1}_{\Int P}\to A^{0,1}_{\frak g}$ by
$$
\bar r_{\ga}((\al_i)) = \ta_1(\al_i) + \ga = \Ad c_i(\al_i) + \ga =
\bar\rh(c_i)(\al_i).\tag1.17
$$
The following statements are easily verified.

\proclaim{Proposition 1.5}
The mapping $\ta$ is an isomorphism of complexes
$(A_{\Ad P},\dif)\to (A_{\frak g},\dif - \ad\ga)$.
\endproclaim

\proclaim{Proposition 1.6}
The following relations hold:
$$
\align
\bar\de_1&\circ\bar r_{\ga} = \ta_2\circ\bar\de_1,\\
\bar r_{\ga}&\circ\bar\rh(g) = \bar\rh(t(g))\circ\bar r_{\ga},\ \
g\in\overline R^0_{\Int P}.
\endalign
$$
Thus, $\bar r_{\ga}$ maps $Z^1(\overline R_{\Int P})$ onto
$Z^1(\overline R_G)$ and induces a bijection $\bar r_{\ga}^*:
H^{0,1}(M,\Int P)\mathbreak\to H^{0,1}(M,G)$ taking $\ep$ to $[\ga]$.
\endproclaim

In what follows, we will consider the case when $E$ is flat, i.e., is
obtained by twisting $M\times G$ with the help of a cocycle $\frak z\in
Z^1(\frak U,\Aut_h G)$. Then the operator $\partial$ is also well-defined in
$\Ph_{\frak e}$, and one can define the subsheaf of antiholomorphic
$q$-forms $\overline\Om^q_{\frak e}\subset\Ph^{0,q}_{\frak e}$. We also
get the bigraded complex $(A_{\frak g},\partial)$ and the Dolbeault
isomorphisms similar to (1.15):
$$
\overline H^{p,q}(M,\frak e)\overset\operatorname{def}\to =
H^p(A_{\frak e}^{*,q},\partial)\simeq H^p(M,\overline\Om^q_{\frak e}).
\tag 1.18
$$

As in the case of the de Rham complex, the correspondence
$G\mapsto\overline R_G$ is a covariant functor from the category of complex
Lie groups into that of cochain complexes. In particular, any holomorphic
homomorphism of complex Lie groups $f: G\to Q$ determines a homomorphism of
sets with distinguished points $f_1^*: H^{0,1}(M,G)\to H^{0,1}(M,Q)$.
\medskip
{\bf 1.3}. There is an important relationship between the de Rham and the
Dolbeault
complexes with values in the same complex Lie group $G$. It corresponds to
the natural inclusion $\iota_G: G\to\Cal O_G$. One checks easily that
the triple $\{\id,\Pi_{0,1},\Pi_{0,2}\}$ is a homomorphism of complexes
$R_G\to\overline R_G$, i.e.,
$$
\align
\bar\de_1&\circ\Pi_{0,1} = \Pi_{0,2}\circ\de_1,\\
\bar\rh(g)&\circ\Pi_{0,1} = \Pi_{0,1}\circ\rh(g),\ \ g\in F_G.
\endalign
$$
Hence we get the following commutative diagram:
$$
\CD
H_{DR}^1(M,G) @>\mu>> H^1(M,G) @>i^*>> H^1(M,\Cal F_G)\\
@V\Pi_{0,1}^*VV @VV\iota_G^*V @|\\
H^{0,1}(M,G) @>>\bar\mu> H^1(M,\Cal O_G) @>>i^*> H^1(M,\Cal F_G),
\endCD\tag1.19
$$
where $\Pi_{0,1}^*$ is induced by $\Pi_{0,1}$.

More generally, suppose that we have a flat group bundle $E$ corresponding
to a cocycle $\frak z\in Z^1(\frak U,\Aut_h G)$. Regarding $E$ as a
holomorphic Lie group bundle, we get the commutative diagram
$$
\CD
H^1_{DR}(M,E) @>\mu>> H^1(M,\Cal C_E) @>i^*>> H^1(M,\Cal E)\\
@V\Pi_{0,1}^*VV @VV\iota_E^*V @|\\
H^{0,1}(M,E) @>>\bar\mu> H^1(M,\Cal O_E) @>>i^*> H^1(M,\Cal E),
\endCD\tag1.20
$$
where $\Pi_{0,1}^*$ is induced by $\Pi_{0,1}$ and $\iota_E^*$ by the natural
inclusion $\iota_E: \Cal C_E\to\Cal O_E$.

We also note that for any $\ga\in Z^1(R_G)$ the following diagram is
commutative:
$$
\CD
A^1_{\Int P} @>r_{\ga}>>  A^1_{\frak g}\\
@V\Pi_{0,1}VV @VV\Pi_{0,1}V \\
A^{0,1}_{\Int P} @>>\bar r_{\ch}> A^{0,1}_{\frak g},
\endCD\tag1.21
$$
where $\Int P$ is the flat group bundle corresponding to $\ga$ and $\ch =
\Pi_{0,1}\ga$.

Now we consider the example of the abelian complex Lie group $G =
\GL_1(\Bbb C) = \Bbb C^{\times}$ which will be useful in the sequel. Note
that $G$ contains the compact subgroup $\U_1 = \{c\in\Bbb C\mid |c| = 1\}$.

\example{Example 1.1}
We want to study the homomorphism of groups $\Pi_{0,1}^*:
H^1_{DR}(M,\Bbb C^{\times})\mathbreak\to H^{0,1}(M,\Bbb C^{\times})$.
Consider the homomorphism $\Exp
: \Bbb C\to\Bbb C^{\times}$ given by $\Exp c =
\exp(2\pi ic)$. Clearly, $\Exp(\Bbb R) = \U_1$. The exact sequences
$$
\align
0&\longrightarrow\Bbb Z\longrightarrow\Bbb R\overset\Exp\to
\longrightarrow\U_1\longrightarrow 1,\\
0&\longrightarrow\Bbb Z\longrightarrow\Bbb C\overset\Exp\to
\longrightarrow\Bbb C^{\times}\longrightarrow 1
\endalign
$$
give rise to the following commutative diagram with exact lines:
$$
\CD
0 @>>> H^1(M,\Bbb Z) @>>> H^1(M,\Bbb R) @>\Exp^*>> H^1(M,\U_1)
@>\de^*>> H^2(M,\Bbb Z)\\
&& @| @VVV @VVV @|\\
0 @>>> H^1(M,\Bbb Z) @>>> H^1(M,\Bbb C) @>\Exp^*>> H^1(M,\Bbb C^{\times})
@>\de^*>> H^2(M,\Bbb Z)\\
&& @| @V\iota_{\Bbb C}^*VV @VV\iota_{(\Bbb C^{\times})^*}V @|\\
0 @>>> H^1(M,\Bbb Z) @>>> H^1(M,\Cal O) @>\Exp^*>> H^1(M,\Cal O^{\times})
@>\de^*>> H^2(M,\Bbb Z).
\endCD
$$
Here $\de^*$ are connecting homomorphisms of the cohomology exact sequences
and arrows without denotation are given by natural inclusions of groups.
It is well known that $\Ker\de^*\subset H^1(M,\Cal O^{\times})$ is precisely
the set of those holomorphic $\GL_1(\Bbb C)$-bundles that are trivial as
smooth principal bundles. By Propositions 1.1 and 1.4, the subgroups
$\Ker\de^*$ may be identified with $H^1_{DR}(M,\U_1),\;
H^1_{DR}(M,\Bbb C^{\times})$ and $H^{0,1}(M,\Bbb C^{\times})$, respectively,
while $\iota_{(\Bbb C^{\times})^*}$ identifies with $\Pi_{0,1}^*$ (see
(1.19)). Using also the classical de Rham and Dolbeault isomorphisms, we
obtain the following commutative diagram with exact lines:
$$
\CD
0 @>>> H^1(M,\Bbb Z) @>>> H^1_{DR}(M,\Bbb R) @>\Exp^*>>
H^1_{DR}(M,\U_1) @>>>\ep\\
&& @| @VVV @VVV \\
0 @>>> H^1(M,\Bbb Z) @>>> H^1_{DR}(M,\Bbb C) @>\Exp^*>>
H^1_{DR}(M,\Bbb C^{\times}) @>>>\ep\\
&& @| @V\Pi_{0,1}^*VV @VV\Pi_{0,1}^* V \\
0 @>>> H^1(M,\Bbb Z) @>>> H^1(M,\Cal O) @>\Exp^*>>
H^{0,1}(M,\Bbb C^{\times}) @>>>\ep.
\endCD
$$
It follows that $H^1_{DR}(M,\Bbb C^{\times})$ and
$H^{0,1}(M,\Bbb C^{\times})$ are connected complex abelian Lie groups,
while $H^1_{DR}(M,\U_1)$ is a torus. The group $\Pic M =
H^{0,1}(M,\Bbb C^{\times})$ is, by definition, the {\it Picard manifold} of
$M$.
\endexample

The following {\it complement problem} seems to be important: under which
conditions is the mapping $\Pi_{0,1}^*: H_{DR}^1(M,G)\to H^{0,1}(M,G)$
surjective? This is the same that for any $\al\in Z^1(\overline R_G)$, to
find a form $\be\in Z^1(R_G)$ such that $\Pi_{0,1}\be =\al$. In
n$^{\circ}2.2$, we will answer this question positively in the case when
$M$ is a compact K\"ahler manifold and $G$ is the Borel subgroup of a
classical complex linear group (e.g., the complex triangular matrix group
$T_n(\Bbb C)$).

\medskip
{\bf 1.4.} Let $\bold V$ be a flat complex vector bundle over a complex
manifold $M$. Then we have the complexes of $\bold V$-valued forms
$(A_{\bold V},d),\;(A_{\bold V}^{*,q},\dif)$ and $(A_{\bold V}^{p,*},\dif)$
(see n$^{\circ}1.2$, where we may formally set $E = \frak e = \bold V$).
Suppose that $M$ is a compact Hermitian manifold and that the structure group
of $\bold V$ is $\U_{n}$. Then we can define a flat Hermitian metric on
$\bold V$. This gives rise to the Laplace operators $\De,\,\De_{\partial}$
and $\De_{\dif}$ in $A_\bold V$ given by
$$
\Delta = dd^* + d^*d,\ \ \De_{\partial} = \partial\partial^*
+\partial^*\partial,\ \ \De_{\dif} = \dif\dif^* + \dif^*\dif,
$$
where $d^*,\,\partial^*,\,\dif^*$ are formally conjugates of $d,\,\partial,\,
\dif$, respectively. Let us denote
$$
\align
\Bbb H^r &= (\Ker\De)\cap A^r_{\bold V},\\
\Bbb H_{\partial}^{p,q} &= (\Ker\De_{\partial})\cap A^{p,q}_{\bold V},\\
\Bbb H_{\dif}^{p,q} &= (\Ker\De_{\dif})\cap A^{p,q}_{\bold V}
\endalign
$$
the vector spaces of harmonic forms. Any harmonic form $\om$ is,
respectively, $d$-, $\partial$- or $\dif$-closed. By Hodge theorem, the
correspondence $\om\mapsto [\om]$ gives the following isomorphisms:
$$
\align
\Bbb H^r &\simeq H^r_{DR}(M,\bold V)\simeq H^r(M,\Cal C_{\bold V}),\\
\Bbb H_{\partial}^{p,q} &\simeq\overline H^{p,q}(M,\bold V)\simeq
\overline H^q(M,\Om^p_{\bold V}),\\
\Bbb H_{\dif}^{p,q} &\simeq H^{p,q}(M,\bold V)\simeq H^q(M,\Om^p_{\bold V}).
\endalign
$$
We shall use the following fact which is well known in the case of the
trivial vector bundle $\bold V$ over a compact K\"ahler manifold (in the
more general case of a harmonic flat vector bundle a proof is outlined
in [ABCKT], p. 104).
 
\proclaim{Proposition 1.7}
Suppose that $M$ is a compact K\"ahler manifold. Then
$$
\Delta = 2\De_{\partial} = 2\De_{\dif}.
$$
It follows that
$$
\Bbb H^r = \bigoplus_{p+q=r}\Bbb H_{\dif}^{p,q},\ \ \Bbb H_{\dif}^{p,q} =
\Bbb H_{\partial}^{p,q}.
$$
\endproclaim
 
In what follows, we will omit the superscripts $\partial$ and $\dif$ in
the notation of harmonic forms.
 
\proclaim{Corollary 1}
Under the above conditions, the vector spaces $\Bbb H^{p,0}$ and
$\Bbb H^{0,q}$ coincide with those of holomorphic $\bold V$-valued $p$-forms
and antiholomorphic $\bold V$-valued $q$-forms, respectively.
\endproclaim
 
\proclaim{Corollary 2}
Under the above conditions,
$$
\align
H^r_{DR}(M,\bold V) &\simeq\bigoplus_{p+q=r} H^{p,q}(M,\bold V)
\simeq\bigoplus_{p+q=r}\overline H^{p,q}(M,\bold V),\\
H^r(M,\Cal C_{\bold V}) &\simeq\bigoplus_{p+q=r}H^p(M,\Om^q_{\bold V})
\simeq\bigoplus_{p+q=r}H^p(M,\overline\Om^q_{\bold V}).
\endalign
$$
\endproclaim
 
Arguing as in the proof of the lemma in [GH], Ch.1, Section 2 (see also
[DGMS]), and applying Proposition 1.7, we obtain the following
`$\partial\dif$-lemma':
 
\proclaim{Corollary 3}
Let $\bold V$ be a flat vector bundle with structure group $\U_{n}$ over a
compact K\"ahler manifold. Suppose that $\om\in A^{p,q}_{\bold V}$ is
$d$-closed and $\partial$- or $\dif$-exact, where $p,q\ge 1$. Then there
exists $\ps\in A^{p-1,q-1}_{\bold V}$ such that $\om =\partial\dif\ps$.
\endproclaim
 
Let $G\subset\GL_n(\Bbb C)$ be a linear complex Lie group and $z = (z_{ij})
\in Z^1(\frak U,G)$. Consider the flat group bundle $E$ corresponding to the
cocycle $\Int z$ (see n$^{\circ}1.1$). Clearly, $z$ determines a flat vector
bundle $\bold V$ with fibre $\Bbb C^n$. Then $E$ can be regarded as a
subbundle of the flat vector bundle $\bold{End}\bold V$ with fibre
$\End\Bbb C^n$ which is obtained by twisting $M\times\End\Bbb C^n$ with the
help of $\Int z$. The group $F_E$ will coincide with a subgroup of the group
$\Aut\bold V$ of automorphisms of $\bold V$.
 
\proclaim{Corollary 4}
Suppose that $z_{ij}\in\U_n$ and that $M$ is a compact K\"ahler manifold. If
$a\in F_E$ satisfies $\Pi_{1,0}\de_0(a) = 0$ or $\Pi_{0,1}\de_0(a) = 0$, then
$\de_0(a) = 0$, i.e., $a$ is a flat section.
\endproclaim
\demo{Proof}
Apply Corollary 1 to the vector bundle $\bold{End}\bold V$.
\enddemo
We return now to the example $G = \Bbb C^{\times}$ (see Section 1.2). Using
Hodge theory (in the case of the trivial bundle $M\times\Bbb C$), we prove
the following
 
\proclaim{Proposition 1.8}
Let $M$ be a compact K\"ahler manifold. The mapping $\Pi_{0,1}^*:
H^1_{DR}(M,\U_1)\to H^{0,1}(M,\Bbb C^{\times})$ is an isomorphism of real
Lie groups.
\endproclaim
\demo{Proof}
Due to the latter diagram of Example 1.1, it suffices to prove that
$\Pi_{0,1}^*: H^1_{DR}(M,\Bbb R)\to H^{0,1}(M,\Bbb C)$ is an isomorphism of
real vector spaces. But this is a
well-known consequence of Proposition 1.7. The inverse to $\Pi_{0,1}^*$ is
given by $\om\mapsto\om + \bar\om,\; \om\in\Bbb H^{0,1}$.
 
\enddemo
 
As a corollary, we get the following well-known fact: for any compact
K\"ahler manifold $M$, the group $\Pic M$ is a complex torus.
 
The above proof shows that any class from $H^1_{DR}(M,\U_1),\;
H^1_{DR}(M,\Bbb C^{\times})$ or \linebreak $H^{0,1}(M,\Bbb C^{\times})$
contains a harmonic form. On the cochain level,
the inverse to $\Pi_{0,1}^*: H^1_{DR}(M,\U_1)\to H^{0,1}(M,\Bbb C^{\times})$
can be given by $\om\mapsto\om - \bar\om,\; \om\in\Bbb H^{0,1}$.
 
\medskip
{\bf 1.5.}
Suppose that we have the semi-direct product $G = B\rtimes A$ of complex Lie
groups. Here $B$ may be regarded as a normal Lie subgroup and $A$ as a Lie
subgroup of $G$. Denote by $p: G\to A$ the natural projection homomorphism
and by $q: A\to G$ the natural inclusion. Then $B = \Ker p$ and $p\circ q =
\id$. Clearly, we have the semi-direct decomposition $\frak g = \frak b
+\!\!\!\!\!\!\supset\frak a$, where $\frak b$ is an ideal and $\frak a$ a
subalgebra of $\frak g$. It follows that $A_{\frak g} = A_{\frak b}
+\!\!\!\!\!\!\supset A_{\frak a}$ for an arbitrary manifold $M$. Any form
$\om\in A^p_{\frak g}$ on $M$ is uniquely decomposed as $\om =\om_{\frak b} +
\om_{\frak a}$, where $\om_{\frak b}\in A^p_{\frak b}$ and
$\om_{\frak a} = dp\circ\om\in A^p_{\frak a}$.
 
The relations between the de Rham (and Dolbeault) cohomology with values in
$G,\,A,\,B$ are described in [O4], Theorems 2.2 and 2.4. In particular,
there is the following commutative diagram with exact lines:
$$
\CD
H_{DR}^1(M,B) @>i_1^*>> H_{DR}^1(M,G) @>p_1^*>>
H_{DR}^1(M,A) @>>> \ep\\
 @V\Pi_{0,1}^*VV @V\Pi_{0,1}^*VV @VV\Pi_{0,1}^*V \\
H^{0,1}(M,B) @>>i_1^*>H^{0,1}(M,G)
@>>p_1^*> H^{0,1}(M,A) @>>> \ep.
\endCD\tag1.22
$$
Both mappings $i_1^*$ are induced by the natural inclusion $i: B\to G$.
We also have the injective mappings $q^*_1: H_{DR}^1(M,A)\to
H_{DR}^1(M,G)$ and $q^*_1: H^{0,1}(M,A)\to H^{0,1}(M,G)$ such that
$p^*_1\circ q^*_1 =\id$. We will identify $H_{DR}^1(M,A)$ and
$H^{0,1}(M,A)$ with their images in $H_{DR}^1(M,G)$ and $H^{0,1}(M,G)$,
respectively.
 
Suppose now that we are given a form $\ga\in Z^1(R_A)$. Denote by $\xi\in
H_{DR}^1(M,A)$ and $\ze\in H^{0,1}(M,A)$ the cohomology classes of $\ga$ and
$\chi = \Pi_{0,1}\ga$. Then $\Pi_{0,1}^*\xi = \ze$. We want to describe the
subsets $(p_1^*)^{-1}(\xi)\subset H_{DR}^1(M,G)$ and
$(p_1^*)^{-1}(\ze)\subset
H^{0,1}(M,G)$ using splitting constructions. As above, choose an open cover
$\frak U = (U_i)$ of $M$ such that $\ga = \de_0(a_i)$ in $U_i$, where $a_i\in
\Ga(U_i,\Cal F_A)$, and consider the 1-cocycle $z = (z_{ij})\in
Z^1(\frak U,A)$ given by $z_{ij} = a_i^{-1}a_j$. Denote by $E$ the flat group
bundle with the $G$ determined by the cocycle $(\Int z_{ij})$. Since
$z_{ij}\in A$, we have the semi-direct decomposition $E = E_B\rtimes E_A$,
where $E_B$ and $E_A$ are group subbundles with fibres $B$ and $A$,
respectively. We also consider the corresponding Lie algebra bundle
$\frak e = \frak e_{\frak b} +\!\!\!\!\!\!\supset \frak e_{\frak a}$.
Then we have the following commutative diagram with exact lines:
$$
\CD
H_{DR}^1(M,E_B) @>i_1^*>> H_{DR}^1(M,E) @>p_1^*>> H_{DR}^1(M,E_A) @>>> \ep\\
 @V\Pi_{0,1}^*VV @V\Pi_{0,1}^*VV @VV\Pi_{0,1}^*V \\
H^{0,1}(M,E_B) @>>i_1^*>H^{0,1}(M,E) @>>p_1^*> H^{0,1}(M,E_A) @>>> \ep.
\endCD\tag1.23
$$
Its right square is got from the right square of (1.22) by applying to the
lines the bijections $r_{\ga}^{-1}$ and $\bar r_{\chi}^{-1}$ (see (1.21)).
Thus, $(p_1^*)^{-1}(\xi)$ is identified with $i_1^*(H_{DR}^1(M,E_B))$, while
$(p_1^*)^{-1}(\ze)$ is identified with $i_1^*(H^{0,1}(M,E_B))$.

\head
2. Cohomology with values in certain solvable algebraic groups
\endhead
 
{\bf 2.1.} Let $G$ be a connected solvable complex linear algebraic group. It
is well known (see, e.g., [Hu], Ch. VII) that $G$ admits the semi-direct
decomposition
$G = N\rtimes S$, where $N$ is a normal unipotent algebraic subgroup (the
unipotent radical) of $G$ and $S\simeq (\Bbb C^{\times})^n$ an algebraic
torus. Respectively, for the tangent Lie algebra $\frak g$ of $G$ we have the
semi-direct decomposition
$\frak g = \frak n +\!\!\!\!\!\!\supset\frak s$, where $\frak n$ and
$\frak s$ are the ideal and the subalgebra of $\frak g$ corresponding to $N$
and $S$, respectively. Denote by $p: G\to S$ the natural projection and by
$q: S\to G$ the natural inclusion.
 
Let $K$ be the compact real form of $S$. Clearly, $K\simeq\U^n_1$.
We will also consider the real algebraic subgroup $G_K = N\rtimes K$ of $G$.
 
The most important example is the subgroup $G = T_n(\Bbb C)$ of
$\GL_n(\Bbb C)$ consisting of upper triangular matrices. In this case
$N = T^0_n(\Bbb C)$ is the normal subgroup of unipotent upper triangular
matrices and $S = D_n(\Bbb C)$ the subgroup of diagonal matrices which is
naturally identified with $(\Bbb C^{\times})^n$. The Lie algebra $\frak g =
\frak t_n(\Bbb C)$ is the subalgebra of upper triangular matrices in
$\frak{gl}_n(\Bbb C)$. Also $\frak n$ is the ideal of nilpotent upper
triangular matrices and $\frak s$ is the subalgebra of diagonal matrices
which is naturally identified with $\Bbb C^n$. Further, $K$ is the subgroup
of unitary diagonal matrices, identified with $\U^n_1$, and its Lie
subalgebra $\frak k\subset\frak s$ is the subalgebra of pure imaginary
diagonal matrices. The de Rham complex $A_{\frak u}$ is identified with the
complex of $\frak s$-valued forms $\ga$ satisfying $\overline{\ga} = -\ga$.
The subgroup $G_K = T_n(\Bbb C)_K$ of $T_n(\Bbb C)$ is formed by matrices
with unitary diagonal.
 
Let $G$ be an arbitrary connected solvable complex linear algebraic group and
suppose that $V$ is an algebraic $G$-module, i.e., that a polynomial
representation $\phi: G\to\GL(V)$ is given. Let $\frak U$ be an open cover of
a complex manifold $M$ and $z\in Z^1(\frak U,G_K)$. Consider the
corresponding flat vector bundle
$\bold V(z)$ of rank $n$ over $M$ with structure group $G_K$. Then there is
the de Rham complex $(A_{\bold V(z)},d)$, as well as the Dolbeault complexes
$(A^{p,*}_{\bold V(z)},\dif)$ and $(A^{*,q}_{\bold V(z)},\partial)$.
 
The projection $p: G_K\to K$ gives rise to the mapping
$p_1^*: H^1(M,G_K)\to H^1(M,K)$. It can be easily described on the cocycle
level. Namely, for
the above cocycle $z = (z_{ij})$ we have $z_{ij} = u_{ij}n_{ij}$, where
$u_{ij}\in K,\; n_{ij}\in N$. Since $u_{ij} = p(z_{ij})$, we see that
$u\in Z^1(\frak U,K)$ and $p_1^*$ sends the cohomology class of $z$ to that
of $u$.
 
Suppose now that  $z_{ij} = g_i^{-1}g_j$ in $U_i\cap U_j\ne\emptyset$, where
$g_i: U_i\to G_K$ are smooth functions. Then $\bold V(z)$ is trivial as a
smooth vector bundle with structure group $G_K$. We have $g_i = u_in_i$ with
$u_i: U_i\to K,\; n_i: U_i\to N$, whence
$$
\aligned
u_{ij} &= p(g_i)^{-1}p(g_j) = u_i^{-1}u_j,\\
z_{ij} &= n_i^{-1}u_{ij}n_j,
\endaligned\tag2.1
$$
In particular, we see that the 0-cochain $(n_i)$ determines an isomorphism
of the flat vector bundles $\bold V(z)$ and $\bold V(u)$ in the category of
smooth vector bundles with structure group $G_K$. In the next lemma, this
cochain should satisfy the condition $\Pi_{0,1}\de_0(n_i) = 0$ for all $i$
which means that the corresponding isomorphism is antiholomorphic.
 
\proclaim{Lemma 2.1}
Let $z\in Z^1(\frak U,G_K)$ be a \v Cech cocycle on a complex
manifold $M$ satisfying the triviality condition described above.
Suppose that $M$ is a compact K\"ahler manifold and that $\Pi_{0,1}\de_0(n_i)
= 0$ for any $i$. Then
 
{\rm (i)} Any class from $H^{0,q}(M,\bold V(z)),\; q\ge 0$, contains
a $d$-closed $\bold V(z)$-valued $(0,q)$-form;
 
{\rm (ii)} Any holomorphic $\bold V(z)$-valued form $\al$ is $d$-closed. If,
in addition, $\al$ is $\partial$-exact, then $\al = 0$;
 
{\rm (iii)} If $\al\in A^{p,0}_{\bold V(z)}$ satisfies $\partial\al = 0$ and
$\dif\al = \partial\be$ for a certain $\be\in A^{p-1,1}_{\bold V(z)}$, then
there exists $\ga\in A^{p-1,0}_{\bold V(z)}$ such that $\al + \partial\ga$ is
holomorphic.
\endproclaim
 
\demo{Proof}
We will use the induction on $n = \dim V$. Assume that $n = 1$. Clearly,
$\phi(G_K)\subset U_1$, and hence $\bold V(z)$ is the flat line bundle
determined by the cocycle $(\phi(z_{ij}))$, where $\phi(z_{ij})\in \U_1$.
Therefore we can apply Hodge theory discussed in Section 1.4. Clearly, (i) is
implied by the fact that any Dolbeault
cohomology class contains a harmonic form (see Proposition 1.7), and
(ii) follows from Corollary 1 of Proposition 1.7. To prove (iii), we note
that $d\partial\be = d\dif\al = 0$. Applying Corollary 3 of the same
proposition, we see that $\partial\be = \dif\partial\ga$, whence
$\dif(\al -\partial\ga) = 0$.
 
To argue by induction, we need the following construction. By Lie theorem,
there exists a basis $e_1,\ldots,e_n$ of $V$ such that $\phi$ is expressed by
the matrix
$$
\phi(g) = \pmatrix \phi_0(g) & *\\ 0 & \tilde\phi(g)\endpmatrix, \ \ g\in G,
$$
where $\tilde\phi: G\to T_{n-1}(\Bbb C)$ is a representation satisfying
$\tilde\phi(g)\in D_{n-1}(\Bbb C)$ for $g\in S$, and hence
$\tilde\phi(g)\in\U_1^{n-1}$, for $g\in K$. The exact sequence of $G$-modules
$$
0\longrightarrow V^0\longrightarrow V\longrightarrow\tilde V
\longrightarrow 0,
$$
where $V^0 = \langle e_1\rangle,\; \tilde V = \langle e_2,\ldots,e_n\rangle$,
gives rise to the following two exact sequences of flat vector bundles:
$$
\align
0&\longrightarrow \bold V^0(z)\overset\la\to\longrightarrow\bold V(z)
\overset\ka\to\longrightarrow \tilde{\bold V}(z)\longrightarrow 0,\tag2.2\\
0&\longrightarrow \bold V^0(z)\overset\la\to\longrightarrow \bold V(u)
\overset\ka\to\longrightarrow\tilde{\bold V}(u)\longrightarrow 0.\tag2.3
\endalign
$$
Since the matrices $\phi(u_{ij})$ are diagonal, the sequence (2.3) is split.
 
From (2.2) we obtain, for any $q\ge 0$, the following exact sequence of
sheaves of antiholomorphic $q$-forms:
$$
0\longrightarrow\overline\Om^q_{\bold V^0(z)}\overset\la\to\longrightarrow
\overline\Om^q_{\bold V(z)}\overset\ka\to\longrightarrow
\overline\Om^q_{\tilde{\bold V}(z)}\longrightarrow 0.\tag2.4
$$
The mappings $n_i$ (see (2.1)) determine an antiholomorphic isomorphism of
(2.2) onto the split exact sequence (2.3). It follows that the exact
sequence of sheaves (2.4) is split. Hence, for
any $p\ge 0$, we get the following split exact sequence of cohomology groups:
$$
0\longrightarrow
H^p(M,\overline\Om^q_{\bold V^0(z)})\overset\la^*\to\longrightarrow
H^p(M,\overline\Om^q_{\bold V(z)})\overset\ka^*\to\longrightarrow
H^p(M,\overline\Om^q_{\tilde{\bold V}(z)})\longrightarrow 0.\tag2.5
$$
Using the isomorphism (1.18), we can interprete (2.5) as the following split
exact sequence:
$$
0\longrightarrow
\overline H^{p,q}(M,\bold V^0(z))\overset\la^*\to\longrightarrow
\overline H^{p,q}(M,\bold V(z))\overset\ka^*\to\longrightarrow
\overline H^{p,q}(M,\tilde{\bold V}(z))\longrightarrow 0.\tag2.6
$$
Let us also note that (2.2) gives rise to the following split exact
sequences:
$$
0\longrightarrow A^{p,q}_{\bold V^0(z)}\overset\la\to\longrightarrow
A^{p,q}_{\bold V(z)}\overset\ka\to\longrightarrow
A^{p,q}_{\tilde{\bold V}(z)}\longrightarrow 0.\tag2.7
$$
 
Now we pass to the induction argument. Suppose that (i), (ii) and (iii) are
true for vector bundles of rank $n-1$.
 
(i) We will assume that $q\ge 1$, since (i) coincides with (ii) in the case
$q = 0$. Let $\al$ be an $\bold V(z)$-valued $(0,q)$-form such that $\dif\al
= 0$, where $\dim V = n$. Consider the exact sequence (2.7) for $p = 0$. By
induction hypothesis, $\ka(\al) = \be +\dif\ph$, where $d\be = 0$ and
$\ph\in A^{0,q-1}_{\tilde{\bold V}(z)}$. Then $\ph =\ka(\ps)$ for a certain
$\ps\in A^{0,q-1}_{\bold V(z)}$. Therefore $\be = \ka(\al -\dif\ps)$. Hence
$\ka(\partial(\al -\dif\ps)) = 0$, and
$\partial(\al -\dif\ps)  = \la(\al')$, where $\al'$ is a $d$-closed
$\bold V^0(z)$-valued $(1,q)$-form. Using (2.6), we conclude that $\al'$ is
 
$\partial$-exact. Applying Corollary 3 of Proposition 1.7, we see that
$\al' = \partial\dif\om$
for a certain $\om\in A^{0,q-1}_{\bold V^0(z)}$. It follows that
$\partial(\al -\dif(\ps+\la(\om))) = 0$. Thus, $\al -\dif(\ps +\la(\om))$ is
the desired $d$-closed form.
 
(ii) Let $\al$ be a holomorphic $\bold V(z)$-valued $p$-form, where
$\dim V = n$. By induction hypothesis, $d\ka(\al) = 0$, whence
$\ka(d\al) = \ka(\partial\al) = 0$. It follows that $\partial\al = \la(\be)$,
where $\be$ is a $\bold V^0(z)$-valued holomorphic $(p+1)$-form. From (2.6)
we
conclude that $\be$ is $\partial$-exact. Using Corollary 1 of Proposition 1.7
applied to $\bold V^0(z)$-valued forms, we see that $\be$ is harmonic, and
hence $\be = 0$. Therefore $\partial\al = d\al = 0$.
 
Suppose that $\al = \partial\ga$. Then $\ka(\al) = \partial\ka(\ga) = 0$ by
induction hypothesis. Therefore $\al = \la(\be)$, where $\be$ is a
$\bold V^0(z)$-valued holomorphic $p$-form. Then (2.6) implies that $\be$ is
$\partial$-exact in $A_{\bold V^0(z)}$, and hence $\be = 0$.
 
(iii) Suppose that $\partial\al = 0$ and $\dif\al = \partial\be$. Clearly,
$\ka(\al)$ satisfies similar conditions. By induction hypothesis,
$\dif(\ka(\al) + \partial\ga_1) = 0$ for a certain $\ga_1\in
A^{p-1,0}_{\tilde{\bold V}(z)}$. By (2.7), we may choose $\ga_2\in
A^{p-1,0}_{\bold V(z)}$
such that $\ka(\ga_2) = \ga_1$. Then $\ka(\dif(\al + \partial\ga_2)) = 0$,
whence $\dif(\al + \partial\ga_2) = \la(\ph)$ for a $\bold V^0(z)$-valued
$(p,1)$-form $\ph$. Clearly, $d\ph = 0$, and (2.6) implies that $\ph$ is
$\partial$-exact. By Corollary 3 of Proposition 1.7, $\ph =
\dif\partial\ga_3$. Then $\dif(\al + \partial(\ga_2 - \ga_3)) = 0$.
\enddemo
 
\proclaim{Corollary}
Let $E$ be the flat group bundle with fibre $G$ determined by the cocycle
$\Int z$. If $a\in F_E$ is holomorphic, then $a$ is flat.
\endproclaim
\demo{Proof}
The argument is similar to the proof of Corollary 4 of Proposition 1.7.
Namely, we fix a faithful holomorphic representation $\phi$ of $G$ and
consider the corresponding flat vector bundle $\bold V(z)$. Then $E$ is a
subbundle of the flat vector bundle $\bold{End}\,\bold V(z)$ determined by
the cocycle $\Int\phi(z)$. Then we apply the statement (ii) to the latter
flat vector bundle.
\enddemo
 
\medskip
{\bf 2.2.} Let $M$ be a compact K\"ahler manifold. Then any class from
$H^{0,1}(M,\Bbb C)$ contains a $d$-closed (i.e., antiholomorphic)
$(0,1)$-form. In this subsection, we consider a non-abelian analogue of this
fact involving a class of solvable algebraic groups. This study is closely
related to the complement problem mentioned in n$^{\circ}1.3$.
 
We retain the notation of n$^{\circ}2.1$. Let $G = N\rtimes S$ be a connected
solvable complex linear algebraic group. Then we have the following
commutative diagram (see (1.22)):
$$
\CD
H^1_{DR}(M,G) @>p^*_1>> H^1_{DR}(M,S)\\
@V\Pi^*_{0,1} VV @VV\Pi^*_{0,1} V \\
H^{0,1}(M,G) @>>p^*_1> H^{0,1}(M,S).
\endCD\tag2.8
$$
As we saw in n$^{\circ}1.5$, both $p^*_1$ are surjective, and $H^1_{DR}(M,S)$
and $H^{0,1}(M,S)$ are naturally embedded in $H^1_{DR}(M,G)$ and
$H^{0,1}(M,G)$, respectively. Since $K$ is a direct factor of $S$, we also
have an embedding $H^1_{DR}(M,K)\subset H^1_{DR}(M,S)$.
 
Due to Proposition 1.8, $\Pi^*_{0,1}$ induces an isomorphism of groups
$H^1_{DR}(M,U)\to H^{0,1}(M,S)$ (which gives the solution of the complement
problem in the case $G = S$). It can be expressed by the isomorphism of the
spaces of $\frak k$-valued and $\frak s$-valued harmonic forms
$\Bbb H^1_{\frak k}\to\Bbb H^{0,1}_{\frak s}$ given by
$-\bar\ch + \ch\mapsto\ch$.
Let us fix a class $\ze\in H^{0,1}(M,S)$ and a harmonic form $\ch\in
A^{0,1}_{\frak s}$ representing this class, and denote $\ga = -\bar\ch +\ch$.
Choose an open cover $\frak U = (U_i)$ of $M$ such that $\ga =
\de_0(u_i)$ in any $U_i$ for certain smooth functions $u_i: U_i\to K$, and
construct the \v Cech cocycle $z = (z_{ij})$, where $z_{ij} = u_i^{-1}u_j$ in
any $U_i\cap U_j\ne\emptyset$. We fix the functions $u_i$, and hence the
cocycle $z$, as well.
 
We will use the twisted complex $R_E$ corresponding
to $\ze$. Here we denote by $E$ the flat group bundle obtained from
$M\times G$ by twisting with the help of $\Int z$, and by $\frak e$ the
corresponding Lie algebra bundle (see n$^{\circ} 1.1$). Since $z$ takes its
values in $K$, we also get the normal group subbundle $E_N\subset E$ with
fibre $N$ corresponding to the cocycle $\Int z|N$ such that $E$ is the
semi-direct product of $N$ and of the trivial bundle $M\times S$.
Respectively, we get the semi-direct decomposition $A_{\frak e} =
A_{\frak e_{\frak n}} +\!\!\!\!\!\!\supset A_{\frak s}$.
 
We also need the mapping $r_{\ga}: A^1_{\frak e}\to A^1_{\frak g}$ introduced
in n$^{\circ} 1.1$. By Proposition 1.3, it induces the bijection $r_{\ga}^*:
H^1_{DR}(M,E)\to H^1_{DR}(M,G)$ such that $r_{\ga}^*(\ep) = \xi$ is the
cohomology class of $\ga$. Similarly (see Proposition 1.6), we have the
mapping $\bar r_{\chi}: A^{0,1}_{\frak e}\to A^{0,1}_{\frak g}$ inducing the
bijection $\bar r_{\chi}^*: H^{0,1}(M,E)\to H^{0,1}(M,G)$ such that
$\bar r_{\chi}^*(\ep) = \ze$.
 
Note that we have the following commutative diagram with exact lines (see
(1.23)):
$$
\CD
H^1_{DR}(M,E_N) @>i^*_1>> H^1_{DR}(M,E) @>p^*_1>> H^1_{DR}(M,S) @>>> \ep\\
 @V\Pi^*_{0,1} VV @V\Pi^*_{0,1} VV @VV\Pi^*_{0,1} V \\
H^{0,1}(M,E_N) @>>i^*_1> H^{0,1}(M,E) @>>p^*_1> H^{0,1}(M,S) @>>> \ep.
\endCD\tag2.9
$$
Here both $p^*_1$ are surjective, and $H^1_{DR}(M,S)$ and $H^{0,1}(M,S)$ are
naturally embedded in $H^1_{DR}(M,E)$ and $H^{0,1}(M,E)$, respectively. As
we saw in n$^{\circ} 1.5$, its right square is got from (1.22) by applying to
the lines the bijections $r_{\ga}^{-1}$ and $\bar r_{\chi}^{-1}$.
 
We will say that $G$ {\it has the Hodge property} if for any compact K\"ahler
manifold $M$ and for any group bundle $E$ over $M$ of the class described
above, the following is satisfied: for each class $\si\in H^{0,1}(M,E)$ such
that $p_1^*(\si) = \ep$ there exists a form $\et\in\si$ with the properties
$$
\et\in A^{0,1}_{\frak e_{\frak n}},\ \ \de_1(\et) = 0.\tag2.10
$$
 
Suppose that $\et$ satisfies (2.10) and denote $\be =
r_{\ga}(\et)\in Z^1(R_G)$. Then $\be = \Ad u_i(\et_i) + \ga$, whence
$$
\be_{\frak s} = \ga,\ \ \Pi_{0,1}\be = \be_{\frak n} + \chi.\tag2.11
$$
Taking into account the remarks made above, we can reformulate our definition
in terms of the complexes $R_G$ and $\overline R_G$. Namely, $G$ has the
Hodge property if for any compact K\"ahler manifold $M$ and for any
$\frak k$-valued harmonic 1-form $\ga$, the following is valid: for each
class $\tilde\si\in H^{0,1}(M,G)$ such that $p_1^*(\tilde\si) =
[\Pi_{0,1}\ga]$, there exists a form $\be\in Z^1(R_G)$ satisfying (2.11).
 
In particular, we see that the Hodge property implies the positive solution
of the complement problem for the group $G$ over a compact K\"ahler manifold.
 
\example{Example 2.1}
Suppose that the unipotent radical $N$ of $G$ is abelian. Then $E$ is a flat
vector bundle with unitary structure group. Hence, Proposition 1.7 implies
that $G$ has the Hodge property.
\endexample
 
Take a form $\et\in R^1_E$ satisfying (2.10). Then the form $\be =
r_{\ga}(\et)$ can be used, instead of $\ga$, for twisting the complex $R_G$.
To do this, we write $\be = \de_0(b_i)$ for a certain smooth function
$b_i: U_i\to G$ in any $U_i$. These functions can be chosen in such a way
that $b_i = u_in_i$ for certain smooth functions $n_i: U_i\to N$. In fact,
$\et_i = \rh(u_i)^{-1}(\be_i) = \rh(u_i)^{-1}(\de_0(b_i)) =
\de_0(u_i^{-1}b_i)$. On the other hand, $\et_i = \de_0(n_i)$ for certain
smooth functions $n_i: U_i\to N$. Therefore $n_i = u_i^{-1}b_ig_i$, where
$g_i: U_i\to G$ satisfy $\de_0(g) = 0$. Then $b_ig_i = u_in_i$ and
$\de_0(b_ig_i) = 0$. Replacing $b_i$ by $b_ig_i$, we get the desired result.
Then we construct the \v Cech cocycle $w = (w_{ij})\in Z^1(\frak U,G)$, where
$$
w_{ij} = b_i^{-1}b_j = n_i^{-1}z_{ij}n_i.
$$
Consider the flat group bundle $E'$ and the flat Lie algebra bundle
$\frak e'$ with fibres $G$ and $\frak g$ determined by the cocycles
$(\Int w_{ij})$ and $(\Ad w_{ij})$, respectively. Clearly, $\be =
\de_0(u_in_i) = \ga + \Ad u_i(\de_0(n_i))$, whence
$$
\et_i = \Ad u_i^{-1}(\be -\ga) = \de_0(n_i),
$$
and $\Pi_{1,0}\de_0(n_i) = 0$. Therefore $\frak e'$ satisfies the conditions
of Lemma 2.1.
 
We will also use the corresponding mapping $r_{\be}: A^1_{\frak e'}\to
A^1_{\frak g}$. Let us find an explicit expression for the mapping
$r_{\be}^{-1}\circ r_{\ga}: \ps\mapsto\ps'$ of the space $A^1_{\frak e}$ onto
$A^1_{\frak e'}$. Clearly, $\rh(u_i)(\ps_i) = r_{\ga}(\ps) = r_{\be}(\ps') =
\rh(u_in_i)(\ps'_i)$, whence
$$
\aligned
\ps_i &= \rh(n_i)\ps'_i,\\
\ps'_i &= \Ad n_i^{-1}(\ps_i - \et_i).
\endaligned\tag2.12
$$
 
Now we are going to prove the Hodge property for certain solvable
algebraic groups. The main role will play the following result.

\proclaim{Lemma 2.2}
Suppose that $G = B\rtimes A$, where $B$ is an abelian unipotent normal
algebraic subgroup and $A$ an algebraic subgroup having the Hodge property.
Then $G$ has the Hodge property.
\endproclaim
 
\demo{Proof}
We have the semi-direct decomposition $A = N_0\rtimes S$, where $N_0$ is the
unipotent radical of $A$ and $S$ is an algebraic torus. Then
$G = N\rtimes S$, where $N = B\rtimes N_0$. Clearly, $N$ is a unipotent
algebraic subgroup, and hence the unipotent radical of $G$.
 
Now let $M$ be a compact K\"ahler manifold. Fix a class $\ze\in H^{0,1}(M,S)$
and a harmonic form $\ch\in A^{0,1}_{\frak s}$ representing this class,
denote $\ga = -\bar\ch +\ch$ and construct the cocycle $z$ and the flat group
bundle $E$ over $M$ as described above. Choose a class $\si\in H^{0,1}(M,E)$
such that $p_1^*(\si) = \ep$. Then $\si = [\al]$ for a certain $\al\in
Z^1(\overline R_{E_N})$. Since $A$ and $B$ are invariant under $\Int z_{ij}$,
we have $E_N = E_B\rtimes E_{N_0}$, where $E_B$ and $E_{N_0}$ are the
subbundles corresponding to the subgroups $B$ and $N_0$ of $N$. Respectively,
$A^{0,1}_{\frak e_{\frak n}} = A^{0,1}_{\frak e_{\frak b}} +
\!\!\!\!\!\!\supset A^{0,1}_{\frak e_{\frak n_0}}$, where
$\frak e_{\frak b}$ and $\frak e_{\frak n_0}$ are the corresponding
subbundles of the Lie algebra bundle $\frak e_{\frak n}$. Hence $\al =
\al_{\frak b} + \al_{\frak e_{\frak n_0}}$, where $\al_{\frak b}\in
A^{0,1}_{\frak e_{\frak b}}$, while $\al_{\frak e_{\frak n_0}}\in
Z^1(\overline R_{E_{N_0}})$. Since $A$ has the Hodge property,
$\al_{\frak e_{\frak n_0}}$ is cohomologous in $\overline R_{E_A}$ to a
form $\et_{\frak a}\in Z^1(R_{E_{N_0}})$. Replacing $\al$ by a cohomologous
cocycle, we may assume that $\al_{\frak e_{\frak n_0}} = \et_{\frak a}$.
 
Consider now the form $\be_{\frak a} = r_{\ga}(\et_{\frak a})\in Z^1(R_A)$.
Clearly, $\et_{\frak a}\in R^1_E$ satisfies (2.10). Hence we may apply the
construction of the flat group bundle $E'$ with fibre $G$ given above using
the form $\be_{\frak a}$ as $\be$. In particular,
we can write $\be_{\frak a} = \de_0(a_i)$ in any $U_i$, where $a_i = u_in_i$
for certain smooth functions $n_i: U_i\to N_0$ satisfying
$\Pi_{1,0}\de_0(n_i) = 0$. Then $E'$ is determined by the cocycle
$w\in Z^1(\frak U,N_0\rtimes K)$ given by $w_{ij} = a_i^{-1}a_j$. Since
$w_{ij}\in A$, we have $E' = E'_B\rtimes E'_A$, where $E'_B$ and $E'_A$ are
the subbundles corresponding to the subgroups $B$ and $A$ of $G$. Using
(2.12), we get
$$
\al' = \Ad n_i^{-1}(\al_i - \et_i) =  \Ad n_i^{-1}((\al_{\frak b})_i),
$$
whence
$$
\bar r^{-1}_{\Pi_{0,1}\be_{\frak a}}(\bar r_{\ch}(\al)) = \Pi_{0,1}\al' =
\al' = \Ad n_i^{-1}(\al_{\frak b})_i).
$$
Thus, $\al'\in Z^1(\overline R_{E'_B})$. Since $B$ is abelian, $E'_B$ is a
flat vector bundle determined by the cocycle $\Int {w_{ij}}|B$ and hence
satisfying the conditions of Lemma 2.1. By the assertion (i) of this lemma,
$\al'$ is cohomologous in $\overline R_{E'_B}$ to a $(0,1)$-form $\al'_1$
such that $d\al'_1 = \de(\al'_1) = 0$. By (2.12), the corresponding form
$\al_1\in Z^1(R_E)$ is expressed by
$$
(\al_1)_i = \rh(n_i)((\al'_1)_i) = \Ad n_i((\al'_1)_i) + \de_0(n_i).
$$
Since $\Pi_{1,0}\de_0(n_i) = 0$, this form satisfies (2.10). Clearly,
$\al_1\in\si$, and Lemma 2.2 is proved.
\enddemo
 
We also will use the following remark.
 
\remark{Remark 2.1}
The definition of the Hodge property immediately implies that certain simple
modifications of a solvable group $G$ preserve the Hodge property.
 
1) Suppose that $G = N\rtimes S$ has the Hodge property.  Then $\tilde G =
N\rtimes\tilde S$ has the Hodge property, too, for any subtorus
$\tilde S\subset S$.
 
2) Let $Z\subset S$ be a central algebraic subgroup of $G = N\rtimes S$. Then
$G$ has the Hodge property if and only if so is $G/Z$.
\endremark
 
\proclaim{Theorem 2.1}
The following connected solvable complex linear algebraic groups $G$ have the
Hodge property:
 
{\rm (i)} $G = T_n(\Bbb C)$,
 
{\rm (ii)} $G$ is the Borel subgroup of a simple complex algebraic group
$\hat G$ of type $A_n,\, B_n,\, C_n,\, D_n,\; E_6$ or $E_7$.
\endproclaim
 
\demo{Proof}
(i) We use the induction on $n$. If $n = 1$, then our assertion is true
by Proposition 1.8. Assume that it is proved for the group $T_{n-1}(\Bbb C)$.
Consider the semi-direct decomposition $G = B\rtimes A$, where
$A\simeq T_{n-1}(\Bbb C)\times\Bbb C^{\times}$ and $B\simeq\Bbb C^{n-1}$
are the subgroups consisting of matrices of the following form:
$$
\align
&A: \pmatrix X & 0\\0 & c\endpmatrix,\ \ X\in T_{n-1}(\Bbb C),\\
&B: \pmatrix I_{n-1} & x\\ 0 & 1 \endpmatrix,\ \ x\in\Bbb C^{n-1}.
\endalign
$$
The induction hypothesis easily implies (see Remark 2.1) that $A$ has the
Hodge property. Hence we may apply Lemma 2.2.
 
(ii) Let $\hat G$ be a simple group of classical type or type $E_6,\, E_7$.
If the type is classical, then we may assume, due to Remark 2.1, that
$\hat G$ is a classical simple group. For the information about parabolic
subgroups used below see, e.g., [He].
 
For $\hat G = \SL_{n+1}(\Bbb C)$, the assertion follows from (i) and Remark
2.1.
 
Suppose that $\hat G = \SO_n(\Bbb C)$. There exists a parabolic subgroup
$P\subset\SO_n(\Bbb C)$ such that
$P = U\rtimes (\SO_{n-2}(\Bbb C)\times\Bbb C^{\times})$, where $U$ is an
abelian unipotent normal subgroup of $P$. It follows that the Borel subgroups
$G_n$ of the groups $\SO_n(\Bbb C)$ satisfy
$G_n = U\rtimes (G_{n-2}\times\Bbb C^{\times})$. Therefore we may argue by
induction on $n$ using Lemma 2.2 and taking into account that $\SO_3(\Bbb C)$
and $\SO_6(\Bbb C)$ are of types $A_1$ and $A_3$, respectively.

Suppose that $\hat G = \Sp_{2n}(\Bbb C)$. There exists a parabolic subgroup
$P\subset\Sp_{2n}(\Bbb C)$ such that $P = U\rtimes\GL_n(\Bbb C)$, where $U$
is an abelian unipotent normal subgroup of $P$. It follows that the Borel
subgroup $G$ of $\Sp_{2n}(\Bbb C)$ satisfies $G = U\rtimes T_n(\Bbb C)$.
Then we apply (i) and Lemma 2.2.
 
Suppose that $\hat G = E_6$. There exists a parabolic subgroup $P\subset E_6$
such that $P = U\rtimes (D_5\times\Bbb C^{\times})$, where $U$
is an abelian unipotent normal subgroup of $P$. It follows that the Borel
subgroup $G$ of $E_6$ satisfies $G = U\rtimes (G_1\times\Bbb C^{\times})$,
where $G_1$ is the Borel subgroup of $D_5$ which has the Hodge
property, as was proved above. Then we apply Lemma 2.2, taking into account
Remark 2.1.
 
Suppose that $\hat G = E_7$. There exists a parabolic subgroup $P\subset E_7$
such that $P = U\rtimes (E_6\times\Bbb C^{\times})$, where $U$ is an abelian
unipotent normal subgroup of $P$. It follows that the Borel subgroup $G$ of
$E_7$ satisfies $G = U\rtimes (G_1\times\Bbb C^{\times})$, where $G_1$ is the
Borel subgroup of $E_6$ which has the Hodge property, as was proved above.
Then we apply Lemma 2.2, taking into account Remark 2.1.
\enddemo
 
\remark{Remark 2.2}
The remaining simple groups $F_4$ and $E_8$ do not contain any parabolic
subgroup with abelian unipotent radical. We do not know, whether the Borel
subgroups of these groups have the Hodge property.
\endremark
 
\medskip

{\bf 2.3.} Here we assume that $M$ is a compact K\"ahler manifold and $G$ is
a connected solvable complex linear algebraic group having the Hodge
property. We will study the cohomology set $H^1_{DR}(M,G)$.
 
We retain the notation introduced in the beginning of n$^{\circ}2.2$. We also
fix a class $\ze\in H^{0,1}(M,S)$ and a harmonic form $\ch\in
A^{0,1}_{\frak s}$ representing this class, and the \v Cech cocycle
$z = (z_{ij})$ with values in $K$ corresponding to the form
$\ga = -\bar\ch +\ch$. To study $H^1_{DR}(M,G)$, it is appropriate to replace
$R_G$ by the twisted complex $R_E$, where $E$ is the flat group bundle
obtained from $M\times G$ by twisting with the help of $\Int z$. The
correspondence between the 1-cohomology of $R_G$ and $R_E$ is given by the
mapping $r_{\ga}: A^1_{\frak e}\to A^1_{\frak g}$ introduced
in n$^{\circ} 1.1$.
 
By definition of the Hodge property, any cohomology class
$\si\in H^{0,1}(M,E)$ such that $p_1^*\si = \ep$ contains a form $\et$
satisfying (2.10).
 
\proclaim{Proposition 2.1}
The form $\et$ is determined uniquely in its cohomology class of $R_E$, up to
a transformation $\Ad a$, where $a$ is a flat section of $E$.
 
The form $\et$ is harmonic, and hence $d\et = [\et,\et] = 0$.
\endproclaim
\demo{Proof}
Let $\et_1\in Z^1(R_{E_N})$ be another $(0,1)$-form such that $\et_1 =
\rh(a)(\et) = \Ad a(\et) + \de_0(a)$, where $a\in F_E$. Then, clearly,
$\Pi_{1,0}\de_0(a) = 0$. By Corollary 4 of Proposition 1.7, $\de_0(a) = 0$.
 
The condition $\de_1(\et) = 0$ implies $\partial\et = 0$. Thus, $\et$ is
antiholomorphic, and hence harmonic due to Corollary 1 of Proposition 1.7. It
follows that $d\et = [\et,\et] = 0$.
\enddemo
 
Now we fix a cohomology class $\si\in H^{0,1}(M,E)$ such that $p_1^*\si =
\ep$ and a form $\et\in\si$ satisfying (2.10).
We are going to study the subset $(\Pi_{0,1}^*)^{-1}(\si)\subset
H^1_{DR}(M,E)$. It is non-empty, containing the cohomology class of $\et$.
A form $\om\in Z^1(R_E)$ will be called {\it canonical} if $\Pi_{0,1}\om\in
Z^1(E_N)$. E.g., $\et$ is canonical.
 
The following important proposition is an easy consequence of the Hodge
property.
 
\proclaim{Proposition 2.2}
Any 1-cohomology class $\ta\in H^1_{DR}(M,E)$ such that $\Pi_{0,1}^*(\ta) =
\si$ contains a cocycle $\om$ that satisfies $\Pi_{0,1}\om = \et$ and hence
is a canonical form.
\endproclaim
\demo{Proof}
By the Hodge property, for any $\al\in\ta$ there exists $a\in F_E$ such that
$$
\et = \bar\rh(a)(\Pi_{0,1}\al) = \Pi_{0,1}\rh(a)(\al) = \Pi_{0,1}\om,
$$
where $\om = \rh(a)(\al)\in\ta$.
\enddemo
 
Consider the form $\be = r_{\ga}(\et)$. In order to get more information on
canonical forms, we will use the twisted complex $R_{E'}$ described in
 
n$^{\circ}2.2$ and the correspondence $r_{\be}^{-1}\circ r_{\ga}:
\ps\mapsto\ps'$ expressed by (2.12).
 
\proclaim{Proposition 2.3}
For any canonical form $\om$ with $\Pi_{0,1}\om = \et$, denote $\al =
\Pi_{1,0}\om$. Then
 
{\rm (i)} $\om'_i = \Ad n_i^{-1}(\al_i)$,
 
{\rm (ii)} $\om'$ is harmonic and satisfies $d\om' = [\om',\om'] = 0$.
 
{\rm (iii)} $d\al = \dif\al = [\et,\al],\; \partial\al = 0$.
\endproclaim
\demo{Proof}
The relation (i) follows easily from (2.12). Since $\om'$ is of type $(1,0)$,
the relation $\de_1(\om')$ implies $\dif\om' = 0$. Thus, $\om'$ is
holomorphic, and hence harmonic by Corollary 1 of Proposition 1.7. This
implies (ii). Now, (iii) follows from (ii) and from the last assertion of
Proposition 2.1.
\enddemo
 
\proclaim{Proposition 2.4}
If $\om$ and $\om_1$ are two canonical forms from the same class $\ta$, then
$\om_1= \Ad f(\om)$ for a certain $f\in H^0(M,\Cal C_E)$.
\endproclaim
\demo{Proof}
By Proposition 2.2, we may assume that $\Pi_{0,1}\om = \Pi_{0,1}\om_1 = \et$.
Suppose that $\om_1= \rh(f)(\om)$, where $f\in F_E$. Using (2.12), we see
that $\om'_1= \rh(f')(\om')$, where $f'\in F_{E'}$ and
$f_i = n_if'_in_i^{-1}$. By Proposition 2.3, $\om'$ and $\om'_1$ are of type
$(1,0)$, and hence $\Pi_{0,1}\de_0(f') = 0$. Due to Corollary of Lemma 2.1,
$\de_0(f') = 0$. Then $\Pi_{1,0}\de_0(f_i) = 0$, whence $\de_0(f_i) = 0$ by
Corollary 4 of Proposition 1.7. Thus, $\om_1= \Ad f(\om)$, where $f$ is a
flat section of $E$.
\enddemo
 
\medskip
 
{\bf 2.4.} We retain the notation introduced above. Here we describe the set
of canonical forms in terms of the complex $(A_{\frak e},d)$.
 
\proclaim{Theorem 2.2}
Any canonical form $\om\in Z^1(R_E)$ can be written as $\om =
\ps + \partial h$, where $h\in A^0_{\frak e}$, $\ps\in A^1_{\frak e}$ is a
uniquely determined harmonic form, and $[\ps,\ps]$ is cohomologous to $0$ in
$(A_{\frak e},d)$.
 
Conversely, let $\ps\in A_{\frak e}$ be a harmonic form such that
$\Pi_{0,1}\ps\in A_{\frak e_{\frak n}}$ and $[\ps,\ps]$ is cohomologous to
$0$ in $(A_{\frak e},d)$. Then $\de_1(\Pi_{1,0}\ps) = \de_1(\Pi_{0,1}\ps) =
0$ and there exists a unique form $\om\in Z^1(R_E)$ such that $\om =
\ps +\partial h$, where $h\in A^0_{\frak e}$. The form $\om$ is canonical.
\endproclaim
 
\demo{Proof} Let $\om$ be a canonical form. We may assume that $\om = \al +
\et$, where $\al = \Pi_{0,1}\om$. By Proposition 2.3 (iii), we have
$\partial\al = 0$. Due to Proposition 1.7, we can write $\al = \ph +
\partial h$, where $h\in A^0_{\frak e}$ and $\ph\in A^{1,0}_{\frak e}$ is
harmonic, i.e., holomorphic. Since $\et$ is harmonic by Proposition 2.1,
$\ps = \ph + \et$ is harmonic, too, and $\om = \ps + \partial h$. The form
$\ph$ is determined uniquely, due to the Hodge decomposition. Clearly,
$$
\align
[\ps,\ps] &= [\om -\partial h,\om -\partial h] =
[\om,\om] + [\partial h,\partial h - 2\om] =
2d\om +\partial[h,\partial h - 2\om]\\
& = \partial(-2\dif h + [h,\partial h - 2\om]).
\endalign
$$
Since $d[\ps,\ps] = 0$, Corollary 3 of Proposition 1.7 implies
that $[\ps,\ps]$ is cohomologous to 0.
 
Conversely, suppose that we have a harmonic form $\ps\in A_{\frak e}$ such
that $\Pi_{0,1}\ps\in A_{\frak e_{\frak n}}$ and $[\ps,\ps] = d\la$, where
$\la\in A^1_{\frak e}$. Let us denote $\ph = \Pi_{1,0}\ps,\; \th =
\Pi_{0,1}\ps$. Evidently,
$$
\align
[\ph,\ph] &= \partial\Pi_{1,0}\la,\tag2.13 \\
[\th,\th] &= \dif\Pi_{0,1}\la,\tag2.14 \\
2[\th,\ph] &= \dif\Pi_{1,0}\la +\partial\Pi_{0,1}\la.\tag2.15
\endalign
$$
Since $[\ph,\ph]$ is holomorphic, Proposition 1.7 and (2.13) imply that
$[\ph,\ph] = \partial\Pi_{1,0}\la = 0$. Similarly, (2.14) implies that
$[\th,\th] = \dif\Pi_{0,1}\la= 0$. In particular, $\de_1(\ph) = \de_1(\th) =
0$. It also follows that $d\dif\Pi_{1,0}\la = d\partial\Pi_{0,1}\la = 0$.
Applying Corollary 3 of Proposition 1.7, we deduce from (2.15) that
$$
[\th,\ph] = \partial\dif g\tag2.16
$$
for a certain $g\in A^0_{\frak e}$.
 
We see, in particular, that $\th$ has the properties of the form $\et$ from
Proposition 2.1. We may denote $\th = \et$ and use the mapping
$r_{\be}^{-1}\circ r_{\ga}: \ps\mapsto\ps'$ described by (2.12).
 
Applying (2.12) to $\ph$ and using (1.3) and (2.16), we get
$$
\Ad n_i(d\ph'_i) = - [\et_i,\ph_i -\et_i] = -\partial\dif g_i,
$$
whence
$$
\aligned
\partial\ph' &= 0,\\
\dif\ph_i &= - \partial(\Ad n_i^{-1}(\dif g_i)),
\endaligned\tag2.17
$$
where $(\Ad n_i^{-1}(\dif g_i))\in A^0_{\frak e'}$.
 
We will seek now for the desired form $\om$. Write $\om = \ps + \partial h =
\al +\et$, where $\al = \ph + \partial h$ and $h\in A^0_{\frak e}$ is to
be chosen in such a way that $\de_1(\om) = 0$. As in Proposition 2.3 (i), we
have $\om'_i = \Ad n_i^{-1}(\al_i)$. Since $\partial\al = 0$, we deduce,
using (1.3), that $\partial\om' = 0$. One can choose $h$ in such a way that
$\dif\om' = 0$. In fact,
$$
\om'_i = \Ad n_i^{-1}(\ph_i) + \Ad n_i^{-1}(\partial h),
$$
where, by (2.12),
$$
\Ad n_i^{-1}(\ph_i) = \ph'_i + \Ad n_i^{-1}(\et_i).
$$
Note that $d(\Ad n_i^{-1}(\et_i)) = 0$. In fact, (2.12) implies that
$\ps' = -\Ad n_i^{-1}(\et_i)$ corresponds to the form $\ps = 0$. Therefore
$\de_1(\ps') = 0$. But $[\ps',\ps'] = 0$ due to Proposition 2.1, whence
$d\ps' = 0$. It follows that $d(\Ad n_i^{-1}(\ph_i)) = d\ph'_i$. From (2.17)
we see that
$\partial(\Ad n_i^{-1}(\ph_i)) = 0$ and $\dif(\Ad n_i^{-1}(\ph_i))$ is
$\partial$-exact. By Lemma 2.1 (iii), there exists $h'\in A^0_{\frak e'}$
such that the form
$$
\Ad n_i^{-1}(\ph_i) + \partial h'_i =
\Ad n_i^{-1}(\ph_i + \partial(\Ad n_i(h'_i))
$$
is holomorphic. Thus, we see that $h_i =\Ad n_i(h'_i)$ satisfy our
conditions.
 
We have chosen $h$ in such a way that $d\om' = 0$. Now, we have
$$
\align
[\om',\om']_i &= \Ad n_i^{-1}([\al_i,\al_i]) =
\Ad n_i^{-1}([\ph_i + \partial h_i,\ph_i + \partial h_i]) \\
&= \partial(\Ad n_i^{-1}([h_i,2\ph_i + \partial h_i])).
\endalign
$$
Since this form is holomorphic, Lemma 2.1 (ii) implies $[\om',\om'] = 0$,
whence $\de_1(\om') = 0$. Therefore $\de_1(\om) = 0$, and hence $\om$ is
canonical.
 
Suppose now that we have another form $\om_1$ satisfying our conditions and
such that $\tilde\om_1 = \ps + \partial h'$ with $h\in A^0_{\frak e}$. It
follows from the above that the form $\partial (h'_i - (h'_1)_i)$, where
$(h'_1)_i = \Ad n_i^{-1}((h_1)_i)$, is holomorphic. By Lemma 2.1 (ii),
$h' = h'_1$, whence $\om = \om_1$.
\enddemo
 
\medskip
 
{\bf 2.5.} Now we are able to formulate our main result. Suppose that $M$ is
a compact K\"ahler manifold. Let us fix a class $\ze\in H^{0,1}(M,S)$, a
harmonic form $\ch$ lying in this class and a \v Cech
cocycle $z = (z_{ij})$ representing the class $\mu([\ga])\in H^1(M,K)$,
where $\ga = -\bar\ch + \ch\in Z^1(R_K)$. We are going to describe the subset
$$
H_{\ze} = (\Pi^*_{0,1})^{-1}((p^*_1)^{-1}(\ze))\subset H^1_{DR}(M,G)
\tag2.18
$$
(see (2.8)). To do this, we note that $r_{\ga}^*$ identifies the right
square of (2.9) with (2.8) (see n$^{\circ}1.5$). In particular, $H_{\ze}$
identifies with
$$
\tilde H_{\ze} = (r_{\ga}^*)^{-1}(H_{\ze}) =
\Ker (p^*_1\circ \Pi^*_{0,1})\subset H^1_{DR}(M,E)
$$
(see (2.9)). To describe the latter subset, we will use harmonic 1-forms from
the de Rham complex $(A_{\frak e},d)$. Note that the vector space
$\Bbb H^p_{\frak e}$ of harmonic $p$-forms is isomorphic to
$H^p_{DR}(M,\frak e)\simeq H^p(M,\Cal C_{\frak e})$ (see n$^{\circ}1.4$).
Since $\Cal C_{\frak e}$ is a sheaf of Lie algebras, there is a natural
bracket operation $[\,,\,]$ in its cohomology, expressed by the bracket of
$\frak e$-valued forms via the de Rham isomorphism. On the other hand, there
is a natural action of the group of flat sections $H^0_{DR}(M,E)\simeq
H^0(M,\Cal C_E)$ on $H^p_{DR}(M,\frak e)\simeq H^p_{DR}(M,\frak e)$ expressed
by the action $\Ad$ on $\frak e$-valued forms and preserving the bracket
operation. Let us define conic algebraic subsets $H^1_{DR}(M,\frak e)_0
\subset H^1_{DR}(M,\frak e)$ and $H^1(M,\Cal C_{\frak e})_0\subset
H^1(M,\Cal C_{\frak e})$ by
$$
\align
H^1_{DR}(M,\frak e)_0 &= \{\xi\in H^1_{DR}(M,\frak e)\mid p_1^*\Pi_{0,1}^*\xi
= 0,\;[\xi,\xi] =0\},\\
H^1(M,\Cal C_{\frak e})_0 &= \{\xi\in H^1(M,\Cal C_{\frak e})\mid
p_1^*\iota_{\frak e}^*(\xi) = 0,\;[\xi,\xi] =0\},
\endalign
$$
where $\iota_{\frak e}: \Cal C_{\frak e}\to\Cal O_{\frak e}$ is the natural
embedding of sheaves (see n$^{\circ}1.3$) and $p_1^*: H^{0,1}(M,\frak e)\to
H^{0,1}(M,\frak s),\; p_1^*: H^1(M,\Cal O_{\frak e})\to
H^1(M,\Cal O_{\frak s})$ are induced by the natural projection $p: \frak e
\to\frak s$. Clearly, they are invariant under $H^0(M,\Cal C_E)$ and
$H^0(M,\Cal C_E)$, respectively. We also need
the corresponding subset of the vector space $\Bbb H^1_{\frak e}$ of
$\frak e$-valued harmonic 1-forms. By Corollary 2 of Proposition 1.7,
$\Bbb H^1_{\frak e} = \Bbb H^{1,0}_{\frak e}\oplus \Bbb H^{0,1}_{\frak e}$,
where the first summand is the set of closed $(1,0)$-forms, while the second
one is the set of closed $(0,1)$-forms. Clearly, $\Bbb H^1_{\frak e}$ is
invariant under $H^0_{DR}(M,E)$. Define the following invariant subset of
this vector space:
$$
\Bbb H^1_{\frak e 0} = \{\ps\in\Bbb H^1_{\frak e}\mid\Pi_{0,1}\ps\in
\Bbb H^{0,1}_{{\frak e}_{\frak n}},\; [\ps,\ps]\text{ is $d$-exact}\}.
$$
 
Note that all these cohomology sets and groups depend, up to isomorphy, on
the class $\ze$ only, being independent of the choice of the cocycles $\ch$
and $z$.
 
\proclaim{Theorem 2.3}
Let $M$ be a compact K\"ahler manifold and $G$ a connected solvable complex
linear algebraic group having the Hodge property. Then
$$
H^1_{DR}(M,G) = \bigsqcup_{\ze\in H^{0,1}(M,S)}H_{\ze},\tag2.19
$$
where $H_{\ze}$ is given by (2.18).
 
Assigning to any cohomology class $\ta\in H_{\ze}$ the harmonic part of
a canonical form that represents the class $(r_{\ga}^*)^{-1}(\ta)\in
\tilde H_{\ze}\subset H^1_{DR}(M,E)$, and to the latter form the
corresponding cohomology class with values in $\Cal C_{\frak e}$, we get the
following bijections:
$$
H_{\ze}\to\Bbb H^1_{\frak e 0}/H^0_{DR}(M,E)\to
H^1_{DR}(M,\frak e)_0/H^0_{DR}(M,E)\to
H^1(M,\Cal C_{\frak e})_0/H^0(M,\Cal C_E).
$$
\endproclaim
\demo{Proof}
The decomposition (2.19) is evident. It is also clear that
$$
\tilde H_{\ze} = \bigsqcup_{\si\in\Ker p_1^*}(\Pi^*_{0,1})^{-1}(\si).
$$
By Proposition 2.2, any cohomology class from $(\Pi^*_{0,1})^{-1}(\si)$
contains a canonical form $\om$. By Theorem 2.2, $\om = \ps + \partial h$,
where $\ps$ is the harmonic part of $\om$ and $h\in A^0_{\frak e}$, and the
correspondence $\om\mapsto\ps$ is a bijection between the sets of canonical
forms from $(\Pi^*_{0,1})^{-1}(\si)$ and forms $\ps\in\Bbb H^1_{\frak e 0}$
satisfying $\Pi_{0,1}\ps\in\si$. By Proposition 2.4, two canonical forms
$\om,\;\om_1$ lie in the same cohomology class if and only if
$\om_1 = \Ad a(\om)$ for a certain $a\in H^0_{DR}(M,E)$. Under this
assumption, $\om_1 = \Ad a(\ps) + \partial(\Ad a(h))$ which implies that
$\ps_1 = \Ad a(\ps)$ is the harmonic part of $\om_1$. Conversely, if
$\om_1 = \ps_1 + \partial h_1$, where $\ps_1 = \Ad a(\ps)$ is harmonic and
$a\in H^0_{DR}(M,E)$, then the canonical forms $\Ad a(\om)$ and $\om_1$ have
the same harmonic part, and hence coincide, due to Theorem 2.2. Thus, we get
a bijection between $\tilde H_{\ze}$ and
$\Bbb H^1_{\frak e 0}/H^0_{DR}(M,E)$. The correspondences between harmonic
forms, de Rham cohomology and cohomology with values in $\Cal C_E$, taking
into account commutativity of (1.20), give bijections between
$\Bbb H^1_{\frak e 0},\; H^1_{DR}(M,\frak e)_0$ and
$H^1(M,\Cal C_{\frak e})_0$ and, hence, between
$\Bbb H^1_{\frak e 0}/H^0_{DR}(M,E),\; H^1_{DR}(M,\frak e)_0/H^0_{DR}(M,E)$
and $H^1(M,\Cal C_{\frak e})_0/H^0(M,\Cal C_E)$.
\enddemo
 
We finish this subsection with two remarks.
 
\remark{Remark 2.3}
The parameter group $H^{0,1}(M,S)$ in (2.19) is isomorphic to the compact
complex torus $(\Pic M)^{\dim S}$ (see Proposition 1.8 and Example 1.1).
\endremark
 
\remark{Remark 2.4}
Let $M$ be a compact K\"ahler manifold, $G$ a complex Lie group, $P$ a
holomorphic principal bundle with base $M$ and fibre $G$. It is known
(see [O1], [O2]) that if $P$ is given by a cocycle $z$ with values in the
maximal compact subgroup $K$ of $G$ (in particular, $P$ is flat), then small
holomorphic deformations of $P$ are flat and are parametrized by
$\Ad P$-valued harmonic (0,1)-forms $\ps$ satisfying $[\ps,\ps] = 0$, where
$\Ad P$ is the Lie algebra bundle determined by the cocycle $\Ad z$. The
condition imposed on $z$ implies that the homomorphism $\pi_1(M)\to G$
corresponding to $P$ takes its values in $K$. If $G$ is linear, we get
a completely reducible representation of $\pi_1(M)$. On the
other hand, it was proved in [GM] that the representation variety
$\Hom(\pi_1(M),\GL_n(\Bbb C))/\Int\GL_n(\Bbb C)$ has at worst quadratic
singularities at the points corresponding to completely reducible
representations. Theorem 2.3 can be regarded as a global result in the same
direction.
\endremark
 
\medskip
 
{\bf 2.6.} In this section, we give some complements to our main result.
 
First, we note that the same methods allow to describe the cohomology set
$H^1_{DR}(M,E_N)$, where $E_N$, as above, denotes the flat group bundle with
fibre $N$ determined by a cocycle $z\in Z^1(\frak U,K)$. To formulate the
result, we introduce algebraic subsets $H^1_{DR}(M,\frak e_{\frak n})_0
\subset H^1_{DR}(M,\frak e_{\frak n}),\; H^1(M,\Cal C_{\frak e_{\frak n}})_0
\subset H^1(M,\Cal C_{\frak e_{\frak n}})$ and $\Bbb H^1_{\frak e_{\frak n}0}
\subset\Bbb H^1_{\frak e_{\frak n}}$ by
$$
\align
H^1_{DR}(M,\frak e_{\frak n})_0 &= \{\xi\in H^1_{DR}(M,\frak e_{\frak n})
\mid [\xi,\xi] = 0\},\\
H^1(M,\Cal C_{\frak e_{\frak n}})_0 &=
\{\xi\in H^1(M,\Cal C_{\frak e_{\frak n}})\mid [\xi,\xi] = 0\},\\
\Bbb H^1_{\frak e_{\frak n}0} &= \{\ps\in\Bbb H^1_{\frak e_{\frak n}}\mid
[\ps,\ps]\text{ is $d$-exact}\}.
\endalign
$$
Now, notice that a class $\ta\in H^1_{DR}(M,E)$ lies in $\Im i_1^* =
\Ker p_1^*$ if and only if any canonical form $\om\in\ta$ satisfies $\om\in
A_{\frak e_{\frak n}}$. In fact, $\om = \om_{\frak n} + \om_{\frak s}$,
where $\om_{\frak n}\in A^1_{\frak e_{\frak n}}$ and $\om_{\frak s}\in
A^{1,0}_{\frak e_{\frak s}}$. If $p_1^*(\ta) = \ep$, then $p(\om) =
\om_{\frak s} = \de_0(a)$, where $a\in F_S$. Clearly, $\Pi_{0,1}\de_0(a) =
0$. If $M$ is a compact, this implies that $\de_0(a) = 0$, and hence $\om =
\om_{\frak n}$. The arguments used to prove Theorem 2.3 yield the following
result:
 
\proclaim{Theorem 2.4}
Let $M$ be a compact K\"ahler manifold, $G$ a connected solvable complex
linear algebraic group having the Hodge property, $N$ its unipotent radical
and $E_N$ the flat group bundle with fibre $N$ determined by a cocycle
$z\in Z^1(\frak U,K)$. Assigning to any cohomology class
$\ta\in H^1_{DR}(M,E_N)$ the harmonic part of a canonical form lying in
$i_1^*(\ta)$, and to the latter form the corresponding de Rham cohomology
class and the cohomology class with
values in $\Cal C_{\frak e_{\frak n}}$, we get the following bijections:
$$
\align
H^1_{DR}(M,E_N)&\to\Bbb H^1_{\frak e_{\frak n}0}/H^0_{DR}(M,E_N)\\
&\to H^1_{DR}(M,\frak e_{\frak n})_0/H^0_{DR}(M,E_N)\to
H^1(M,\Cal C_{\frak e_{\frak n}})_0/H^0(M,\Cal C_{E_N}).
\endalign
$$
\endproclaim
 
\remark{Remark 2.5}
In particular, in the simplest case $z = e$, we get a description of
$H^1_{DR}(M,N)$ in terms of $H^1(M,\frak n)$. Since $N$ is contractible, we
have a bijection between $H^1_{DR}(M,N)$ and $\Hom(\pi_1(M),N)/\Int N$. A
similar description of the latter set for any simply connected complex
nilpotent Lie group $N$ in terms of harmonic forms can be deduced from the
results of [DGMS].
\endremark
 
Next, we give a generalization of Theorem 2.3 based on the fact that the
Hodge property deals with the $(0,1)$-parts of our matrix forms only. This
suggests to consider the following non-abelian cochain complex. Let $\hat G$
be an arbitrary complex linear algebraic group and $G$ its connected solvable
algebraic subgroup. Denote by $\hat\frak g\supset\frak g$ the corresponding
Lie algebras. Consider the triple $R_{\hat G,G} =
\{R_{\hat G,G}^0,R_{\hat G,G}^1,R_{\hat G,G}^2\}$ defined by
$$
R_{\hat G,G}^0 = F_G,\; R_{\hat G,G}^p =
\{\al\in A^p_{\hat\frak g}\mid\Pi_{0,p}\al\in A^{0,p}_{\frak g}\},\; p = 1,2.
$$
One verifies easily that $R_{\hat G,G}$ is a subcomplex (in the sense of
[On3, On4]) of the de Rham complex $R_{\hat G}$ with values in $\hat G$ (see
n$^{\circ}1.1$). This means that $\de_p(R_{\hat G,G}^p)\subset
R_{\hat G,G}^{p+1},\; p = 0,1$, and that $R_{\hat G,G}^p,\; p = 1,2$, are
invariant under the actions $\Ad$ and $\rh$ of $R_{\hat G,G}^0$. On the other
hand, $R_G$ is a subcomplex of $R_{\hat G,G}$ and coincides with
$R_{\hat G,G}$ in the case $\hat G = G$. One defines the 0-cohomology group
$H^0(R_{\hat G,G})$ (which coincides with $H^0(R_G)$) and the 1-cohomology
set $H^1(R_{\hat G,G}) = Z^1(R_{\hat G,G})/\rh(F_G)$.
 
As in n$^{\circ}1.3$, we see that the triple $\{\id,\Pi_{0,1},\Pi_{0,2}\}$ is
a homomorphism of complexes $R_{\hat G,G}\to\overline R_G$. It gives rise to
the homomorphism of sets with distinguished points $\Pi_{0,1}^*:
H^1(R_{\hat G,G})\to H^{0,1}(M,G)$. Suppose that $M$ is a compact K\"ahler
manifold and $G$ has the Hodge property. Then $\Pi_{0,1}^*$ is a surjection.
 
To describe $H^1(R_{\hat G,G})$, we will use a twisting of the complex
$R_{\hat G,G}$. As above,
let us fix a class $\ze\in H^{0,1}(M,S)$, a harmonic form $\ch$ lying in this
class and a \v Cech cocycle $z = (z_{ij})$ representing the class
$\mu([\ga])\in H^1(M,K)$, where $\ga = -\bar\ch + \ch\in Z^1(R_K)$. Denote by
$\hat E$ the flat group bundle obtained from $M\times\hat G$ by twisting with
the help of the cocycle $\Int z$ with values in $\Aut\hat G$, and by
$\hat\frak e$ the corresponding Lie algebra bundle. Clearly, the flat bundles
$E$ and $\frak e$ studied above are subbundles of $\hat E$ and $\hat\frak e$,
respectively. Then we can define the subcomplex $R_{\hat E,E}$ of
$R_{\hat E}$ given by
$$
R_{\hat E,E}^0 = F_E,\; R_{\hat E,E}^p =
\{\al\in A^p_{\hat\frak e}\mid\Pi_{0,p}\al\in A^{0,p}_{\frak e}\},\; p = 1,2.
$$
The triple $\{\id,\Pi_{0,1},\Pi_{0,2}\}$ is a homomorphism of complexes
$R_{\hat E,E}\to\overline R_E$ giving rise to the homomorphism $\Pi_{0,1}^*:
H^1(R_{\hat E,E})\to H^{0,1}(M,E)$. Clearly, we have the following
commutative diagram:
$$
\CD
H^1(R_{\hat E,E}) @>\Pi^*_{0,1}>> H^{0,1}(M,E) @>p^*_1>> H^{0,1}(M,S)\\
@Vr^*_{\ga}VV @V\bar r^*_{\ch}VV @VV\bar r^*_{\ch}V \\
H^1(R_{\hat G,G}) @>>\Pi^*_{0,1}> H^{0,1}(M,G) @>>p^*_1> H^{0,1}(M,S),
\endCD
$$
where $r^*_{\ga}$ is the bijection determined by the mapping $r_{\ga}:
R_{\hat E,E}^1\to R_{\hat G,G}^1$.
 
Define the subset
$$
\hat H_{\ze} = (\Pi^*_{0,1})^{-1}((p^*_1)^{-1}(\ze))\subset
H^1_{DR}(M,\hat G).\tag2.20
$$
Then
$$
(r^*_{\ga})^{-1}(\hat H_{\ze}) = \Ker(p^*_1\circ\Pi^*_{0,1})
\subset H^1_{DR}(M,\hat E).
$$
Choose a $\ta\in (r^*_{\ga})^{-1}(\hat H_{\ze})$. A cocycle $\om\in\ta$ is
called {\it canonical} if $\Pi_{0,1}\om\in Z^1(E_N)$. Our assumptions imply
that any
such class $\ta$ contains a canonical form which is unique up to a
transformation $\Ad g$, where $g$ is a flat section of $E$. Finally, define
the following conic subset of the vector space $\Bbb H^1_{\hat\frak e}$:
 
$$
\Bbb H^1_{\hat\frak e 0} = \{\ps\in\Bbb H^1_{\hat\frak e}\mid\Pi_{0,1}\ps\in
\Bbb H^{0,1}_{{\frak e}_{\frak n}},\; [\ps,\ps]\text{ is $d$-exact}\}.
$$
Then we get the following result:
 
\proclaim{Theorem 2.5}
Let $M$ be a compact K\"ahler manifold and suppose that $G$ has the Hodge
property. Then
$$
H^1(R_{\hat G,G}) = \bigsqcup_{\ze\in H^{0,1}(M,S)}\hat H_{\ze},
$$
where $\hat H_{\ze}$ is given by (2.20).
 
Assigning to any cohomology class $\ta\in\hat H_{\ze}$ the harmonic part of
a canonical form that represents the class $(r_{\ga}^*)^{-1}(\ta)$, we get
the following bijection:
$$
\hat H_{\ze}\to\Bbb H^1_{\hat\frak e 0}/H^0_{DR}(M,E).
$$
\endproclaim
 
The proof goes along the same lines as that of Theorem 2.3.
 
\vskip1cm
\centerline{REFERENCES}
\bigskip
[ABCKT] J. Amor\'os, M. Burger, K. Corlette, D. Kotschick and D. Toledo,
Fundamental Groups of Compact K\"{a}hler Manifolds. Amer. Math. Soc.,
Princeton, N.J. (1996).
 
[Br1] A. Brudnyi, $\dif$-equations on compact K\"{a}hler manifolds
and representations of fundamental groups. Research Thesis. Technion,
Haifa (1995).
 
[Br2] A. Brudnyi, Representations of fundamental groups of compact K\"{a}hler
manifolds in solvable matrix groups. To appear in Michigan Math. J. {\bf 46}
(1999).
 
[DGMS] P. Deligne, P. Griffiths, J. Morgan and D. Sullivan, Real homotopy
theory of K\"{a}hler manifolds. Invent. Math. {\bf 29} (1975), 245--274.
 
[GM] W.M. Goldman and J.J. Millson, The deformation theory of representations
of fundamental groups of compact K\"{a}hler manifolds. Publ. Math. IHES
{\bf 67} (1988), 43--96.
 
[GH] P. Griffiths and J. Harris, Principles of algebraic geometry. J. Wiley
and Sons, N.Y. (1978).
 
[He] S. Helgason, Differential Geometry, Lie Groups, and Symmetric Spaces.
Academic Press, N.Y. e.a. (1978).
 
[Hu] J.E. Humphreys, Linear Algebraic Groups. Springer-Verlag, N.Y. e.a.
(1975).
 
[O1] A.L. Onishchik, Deformations of fiber bundles. DAN  {\bf 161} (1965),
45--47 (in Russian). English transl.:  Doklady 1965, Tom 161, 369--371.
 
[O2] A.L. Onishchik, On deformations of holomorphic fibre bundles. In: Modern
Problems of Analytic Functions Theory. Nauka, Moscow (1966), 236--239
(in Russian).
 
[O3] A.L.Onishchik, Certain concepts and applications of non-abelian
cohomology theory. Transact. of the Moscow Math.Soc. {\bf 17} (1967), 45--88
(in Russian). English transl.: Transact. of the Moscow Math.Soc. {\bf 17},
Amer. Math. Soc. (1969).
 
[O4] A.L.Onishchik, On non-abelian cochain complexes, In: Topics in the Group
Theory and Homological algebra. Yaroslavl University, Yaroslavl (1998), 171--
197.
 
\enddocument